\documentclass[12pt,reqno]{amsart}

\usepackage{amssymb,amsmath}
\usepackage{amsfonts}
\usepackage{amsthm}
\usepackage{mathrsfs,graphicx}
\usepackage{bm}

\theoremstyle{plain} 
\newtheorem{theorem}{Theorem}
\newtheorem{corollary}[theorem]{Corollary}

\theoremstyle{definition} 
\newtheorem{definition}[theorem]{Definition}
\newtheorem{remark}[theorem]{Remark}
\newtheorem{example}[theorem]{Example}
\newcommand{\R}{\ensuremath{\mathbb{R}}}
\newcommand{\Z}{\ensuremath{\mathbb{Z}}}

\newcommand{\N}{\ensuremath{\mathbb{N}}}
\newcommand{\C}{\ensuremath{\mathbb{C}}}

\numberwithin{equation}{section}
\numberwithin{theorem}{section}

\small\normalsize

\begin{document}

\title[Best Ulam constants for two-dimensional linear systems]{Best Ulam constants for two-dimensional nonautonomous linear differential systems}
\author[Anderson]{Douglas R. Anderson} 
\address{Department of Mathematics, 
         Concordia College, 
         Moorhead, MN 56562 USA}
\email{andersod@cord.edu}
\author[Onitsuka]{Masakazu Onitsuka}
\address{Department of Applied Mathematics, 
         Okayama University of Science, 
         Okayama, 700-0005, Japan}
\email{onitsuka@ous.ac.jp}
\author[O'Regan]{Donal O'Regan}
\address{School of Mathematical and Statistical Sciences, 
         University of Galway, 
         Galway, Ireland}
\email{donal.oregan@nuigalway.ie}


\begin{abstract}
This study deals with the Ulam stability of nonautonomous linear differential systems without assuming the condition that they admit an exponential dichotomy. In particular, the best (minimal) Ulam constants for two-dimensional nonautonomous linear differential systems with generalized Jordan normal forms are derived. The obtained results are applicable not only to systems with solutions that exist globally on $(-\infty,\infty)$, but also to systems with solutions that blow up in finite time. New results are included even for constant coefficients. A wealth of examples are presented, and approximations of node, saddle, and focus are proposed. In addition, this is the first study to derive the best Ulam constants for nonautonomous systems other than periodic systems. 
\end{abstract}

\maketitle\thispagestyle{empty}


\section{Introduction}

Let $I$ be an interval of $\R$, and let $n\in \N$. We consider the nonautonomous inhomogeneous linear differential equation
\begin{equation}
 \bm{x}' = A(t)\bm{x}+\bm{f}(t),
 \label{perturbed}
\end{equation}
where $A: I \to \C^{n\times n}$ is a continuous matrix function and $\bm{f}: I \to \C^n$ is a continuous vector function, and $\bm{x}\in \C^n$. When $\bm{f}(t) \equiv \bm{0}$, equation \eqref{perturbed} is reduced to the homogeneous linear differential equation
\begin{equation}
 \bm{x}' = A(t)\bm{x}.
 \label{unperturbed}
\end{equation}
The concept of Ulam stability was proposed by S. M. Ulam in 1940 in the field of functional equations, and introduced to ordinary differential equations in the 1990s (see \cite{AgaXuZha,BrzPopRasXu,BrzPopRasXu2}). Ulam stability for various types of differential equations has been considered in recent years. However, work on vector equations is still in progress, while work on scalar equations with constant coefficients has been extensively explored. We can cite references \cite{BacDra1,BacDra2,BacDraPitSin,BusLupO'R,BusO'RSaiTab,Dragicevic1} as previous studies on \eqref{perturbed} and \eqref{unperturbed}, and nonlinear equations containing them or their discrete versions. It should be noted that these previous studies clarify the relationship between Ulam stability and the property that the linear part admits an exponential dichotomy. The exponential dichotomy is a very useful property that works well with constant coefficient systems and periodic systems. In addition, the relationship between Ulam stability and the shadowing property has recently been clarified. For shadowing properties see \cite{BacDra1,BacDra2,BacDra3,BacDraPitSin,DraPit}. In particular, we find that Ulam stability is a special case of the Lipschitz shadowing property. The definition of Ulam stability used in this study is as follows. 


\begin{definition}\label{Ulamdef}
Equation \eqref{perturbed} is Ulam stable on $I$, if there exists a constant $K>0$ with the following property:
\begin{quote}
For any $\varepsilon>0$, and for any continuously differentiable function $\bm{\phi}: I \to \C^n$ satisfying
\[ \sup_{t \in I}\|\bm{\phi}' - A(t)\bm{\phi} - \bm{f}(t)\| \le \varepsilon, \]
there exists a solution $\bm{x}: I \to \C^n$ of \eqref{perturbed} such that
\[ \sup_{t \in I}\|\bm{\phi}(t)-\bm{x}(t)\| \le K \varepsilon, \]
\end{quote}
where $\|\bm{v}\|$ is the norm of $\bm{v} \in \C^n$. We call such a $K$ an Ulam constant for \eqref{perturbed} on $I$. 
\end{definition}

Let $R: I \to \C^{n\times n}$ be a continuously differentiable regular matrix function; that is, the inverse matrix $R^{-1}(t)$ exists for all $t \in I$. Using $\bm{y}:=R(t)\bm{x}$, equation \eqref{unperturbed} is transformed into the linear differential equation
\begin{equation}
 \bm{y}' = J(t)\bm{y}, \quad J(t) := \left(R'(t)+R(t)A(t)\right)R^{-1}(t).
 \label{Jordan}
\end{equation}
If $R(t)$ is a constant matrix $R$, then $J(t)=RA(t)R^{-1}$, and it is well known that this transformation preserves many qualitative properties between \eqref{unperturbed} and \eqref{Jordan}. For example, in this case it is clear that the zero solution of \eqref{unperturbed} is asymptotically stable if and only if the zero solution of \eqref{Jordan} is asymptotically stable. So, the first question in this study is the following: ``Is Ulam stability preserved in \eqref{unperturbed} and \eqref{Jordan} under appropriate assumptions for $R(t)$?''. This study provides a positive answer to this question. In fact, it can be shown that \eqref{unperturbed} is Ulam stable on $I$ if and only if \eqref{Jordan} is Ulam stable on $I$ when $\|R(t)\|$ and $\|R^{-1}(t)\|$ are bounded (see Theorem \ref{iff2}), where $\|M\|$ denotes the matrix norm of $M \in \C^{n\times n}$. Furthermore, we can see that \eqref{perturbed} is Ulam stable on $I$ if and only if \eqref{unperturbed} is Ulam stable on $I$ without additional conditions (see Theorem \ref{iff1}). These facts mean that the Ulam stability of \eqref{perturbed} reduces to that of \eqref{Jordan}. If we assume that $A(t)$ and $R(t)$ are real-valued constant matrices $A$, $R \in \R^{n\times n}$ with $\det A \ne 0$, and $n=2$, then the matrix $J=RAR^{-1}$ will be one of the following matrices:
\[ \text{(i)}\;\; \begin{pmatrix}
\lambda_1 & 0\\
0 & \lambda_2
\end{pmatrix}, \quad 
\text{(ii)}\;\; \begin{pmatrix}
\lambda & 1\\
0 & \lambda
\end{pmatrix}, \quad 
\text{(iii)}\;\; \begin{pmatrix}
\alpha & \beta\\
-\beta & \alpha
\end{pmatrix}, \]
by choosing an appropriate regular matrix $R$, where $\lambda_1$, $\lambda_2$, $\lambda$, $\beta \in \R\setminus \{0\}$, $\alpha \in \R$. These matrices are known as Jordan normal forms. In particular, case (iii) is called a real Jordan normal form. The following geometric classifications of the origin $\bm{y}=\bm{0}$ of \eqref{Jordan} are well known: 
\begin{itemize}
  \item if $\lambda_1 \ge \lambda_2 > 0$, then the origin of \eqref{Jordan} with matrix (i) is an unstable node;
  \item if $0 > \lambda_1 \ge \lambda_2$, then the origin of \eqref{Jordan} with matrix (i) is a stable node;
  \item if $\lambda_1 > 0 > \lambda_2$, then the origin of \eqref{Jordan} with matrix (i) is a saddle point;
  \item if $\lambda > 0$, then the origin of \eqref{Jordan} with matrix (ii) is an unstable node;
  \item if $0 > \lambda$, then the origin of \eqref{Jordan} with matrix (ii) is a stable node;
  \item if $\alpha > 0$, then the origin of \eqref{Jordan} with matrix (iii) is an unstable focus;
  \item if $0 > \alpha$, then the origin of \eqref{Jordan} with matrix (iii) is a stable focus;
  \item if $\alpha = 0$, then the origin of \eqref{Jordan} with matrix (iii) is a center.
\end{itemize}
As can be seen from Definition \ref{Ulamdef}, to guarantee the Ulam stability of \eqref{Jordan} is to give approximations of node, saddle, and focus points with accuracy. In fact, using the results of \cite{AndOni5}, it is possible to determine whether cases (i)--(iii) are Ulam stable or not. In conclusion, under the assumption of $I=\R$, it is unstable if $\alpha=0$, but Ulam stable otherwise (see, \cite[Theorem 3.1--3.3, and Lemma 4.2.]{AndOni5}). Moreover, if $\det A = 0$, then $A$ has at least one 0 eigenvalue. For this case, it is known that \eqref{unperturbed} is not Ulam stable on $\R$ (see \cite[Lemma 4.1]{AndOni5}.). Therefore, it is possible to determine the Ulam stability of \eqref{perturbed} with real-valued constant matrices in all cases. As mentioned above, Ulam stability is an approximation concept. Thus, the next important question is ``How close are the approximate and true solutions?''. If the minimum Ulam constant exists, it is called the best Ulam constant. 
In recent years, the existence of the best Ulam constant has been known in some constant coefficients linear differential equations, linear difference equations, and functional equations (see, \cite{AndOni3,AndOni5,AndOniRas,BaiBlaPop2,BaiPop1,BaiPop2,PopRas1,PopRas2,PopRas3}). 
Some of the best Ulam constants for \eqref{unperturbed} have been obtained in \cite{AndOni5}, when $A(t)$ is a real-valued constant matrix $A \in \R^{2\times 2}$. In particular, in the case of matrices (i) and (ii), we can obtain the best Ulam constants. Thus, some best Ulam constants for constant coefficients equations have been obtained, whereas, to the best of the authors' knowledge, the existence of best Ulam constants for variable coefficients equations is limited to scalar and periodic coefficient equations (see, \cite{AndOniRas,FukuOni1,FukuOni2}). 
In this study, we focus on this point of view to establish the Ulam stability of Jordan normal forms generalized to variable coefficients, and we aim to derive their best Ulam constants. 
Specifically, the main purpose of this paper is to derive the best Ulam constants of \eqref{unperturbed} whose coefficients are the generalized Jordan normal forms: 
\[ \text{(I)}\;\; \begin{pmatrix}
\lambda_1(t) & 0\\
0 & \lambda_2(t)
\end{pmatrix}, \quad 
\text{(II)}\;\; \begin{pmatrix}
\lambda(t) & \mu(t)\\
0 & \lambda(t)
\end{pmatrix}, \quad 
\text{(III)}\;\; \begin{pmatrix}
\alpha(t) & \beta(t)\\
-\beta(t) & \alpha(t)
\end{pmatrix}, \]
where $\lambda_1$, $\lambda_2$, $\lambda$, $\mu: I \to \C$, and $\alpha$, $\beta: I \to \R$ are continuous functions. In addition, if we restrict these functions to real-valued functions, we propose that Ulam stability ensures that phase portraits of the orbits of the approximate solutions on the phase plane give approximations of node, saddle, and focus points, respectively.


\section{General theory of Ulam stability}

In this section, we show that equivalence relations are obtained for the Ulam stability of the three equations \eqref{perturbed}, \eqref{unperturbed}, and \eqref{Jordan}. The first result is as follows. 

\begin{theorem}\label{iff1}
Equation \eqref{perturbed} is Ulam stable on $I$, with an Ulam constant $K$ if and only if equation \eqref{unperturbed} is Ulam stable on $I$, and an Ulam constant is the same $K$.
\end{theorem}

\begin{proof}
We only have to show that if \eqref{unperturbed} is Ulam stable on $I$, with an Ulam constant $K$, then \eqref{perturbed} is Ulam stable on $I$, with an Ulam constant $K$, since the converse is trivial. Let $\varepsilon>0$. We consider the continuously differentiable function $\bm{\phi}: I \to \C^n$ satisfying
\[ \sup_{t \in I}\|\bm{\phi}' - A(t)\bm{\phi} - \bm{f}(t)\| \le \varepsilon. \]
If $\bm{z}: I \to \C^n$ is a solution of \eqref{perturbed}, then we obtain
\begin{eqnarray*}
	\varepsilon &\ge& \sup_{t \in I}\|\bm{\phi}' - A(t)\bm{\phi} - \bm{f}(t)\|
	= \sup_{t \in I}\|\bm{\phi}' - A(t)\bm{\phi} - (\bm{z}' - A(t)\bm{z})\| \\
	&=& \sup_{t \in I}\|(\bm{\phi}-\bm{z})' - A(t)(\bm{\phi}-\bm{z})\|.
\end{eqnarray*}
This, together with the Ulam stability of \eqref{unperturbed}, yields that there is a solution $\bm{y}: I \to \C^n$ of \eqref{unperturbed} such that
\[ \sup_{t \in I}\|(\bm{\phi}(t)-\bm{z}(t))-\bm{y}(t)\| \le K \varepsilon. \]
Let $\bm{x}:= \bm{y}+\bm{z}$. Then $\bm{x}$ is a solution of \eqref{perturbed}. From
\[ \sup_{t \in I}\|\bm{\phi}(t)-\bm{x}(t)\| = \sup_{t \in I}\|(\bm{\phi}(t)-\bm{z}(t))-\bm{y}(t)\| \le K \varepsilon, \]
we conclude that \eqref{perturbed} is Ulam stable on $I$, with an Ulam constant $K$. This completes the proof.
\end{proof}

The second result is as follows. 

\begin{theorem}\label{uimpliesJ}
Suppose that $\|R(t)\|$ and $\|R^{-1}(t)\|$ are bounded on $I$. If equation \eqref{unperturbed} is Ulam stable on $I$, with an Ulam constant $K$, then equation \eqref{Jordan} is Ulam stable on $I$, and an Ulam constant is $\left(\sup_{t \in I}\|R(t)\|\right)\left(\sup_{t \in I}\|R^{-1}(t)\|\right)K$. 
\end{theorem}

\begin{proof}
Assume that $\|R(t)\|$ and $\|R^{-1}(t)\|$ are bounded on $I$, and equation \eqref{unperturbed} is Ulam stable on $I$, with an Ulam constant $K$. 
Let $\varepsilon >0$, and let $\bm{\eta}: I \to \C^n$ satisfy
$$ \sup_{t \in I}\|\bm{\eta}'-J(t)\bm{\eta}\| \le \varepsilon, $$
where $J(t) = \left(R'(t)+R(t)A(t)\right)R^{-1}(t)$. The purpose of this proof is to show that there exists a solution $\bm{y}: I \to \C^n$ to equation \eqref{Jordan} that satisfies the inequality
$$ \sup_{t \in I}\|\bm{\eta}(t)-\bm{y}(t)\| \le \left(\sup_{t \in I}\|R(t)\|\right)\left(\sup_{t \in I}\|R^{-1}(t)\|\right)K\varepsilon. $$

Let $\bm{\phi} = R^{-1}(t)\bm{\eta}$. Then, by $J(t) = \left(R'(t)+R(t)A(t)\right)R^{-1}(t)$ and $R^{-1}(t)R(t)=I$, we have
$$ A(t) = R^{-1}(t)J(t)R(t)-R^{-1}(t)R'(t) $$
and
$$ (R^{-1}(t))' = -R^{-1}(t)R'(t)R^{-1}(t). $$
Using these equalities, we obtain
\begin{align*}
\sup_{t \in I}\|\bm{\phi}'-A(t)\bm{\phi}\|
  &= \sup_{t \in I}\|R^{-1}(t)\bm{\eta}'+(R^{-1}(t))'\bm{\eta}-A(t)R^{-1}(t)\bm{\eta}\| \\
  &= \sup_{t \in I}\|R^{-1}(t)(\bm{\eta}'-J(t)\bm{\eta})\| \\
  &\le \sup_{t \in I}\|R^{-1}(t)\|\sup_{t \in I}\|\bm{\eta}'-J(t)\bm{\eta}\| \le \varepsilon\sup_{t \in I}\|R^{-1}(t)\|.
\end{align*}
Since $\|R^{-1}(t)\|$ is bounded on $I$, and \eqref{unperturbed} is Ulam stable on $I$, with an Ulam constant $K$, we see that there is a solution $\bm{x}: I \to \C^n$ of \eqref{unperturbed} such that
\begin{equation}
 \sup_{t \in I}\|\bm{\phi}(t)-\bm{x}(t)\| \le K \varepsilon\sup_{t \in I}\|R^{-1}(t)\| < \infty.
 \label{eq:2.02}
\end{equation}
Letting $\bm{y} = R(t)\bm{x}$, we have
$$ \bm{y}' = R(t)\bm{x}'+R'(t)\bm{x} = (R(t)A(t)+R'(t))\bm{x} = J(t)\bm{y}, $$
and so that $\bm{y}$ is a solution to \eqref{Jordan}. In addition, by \eqref{eq:2.02} and the boundedness of $\|R(t)\|$, we obtain
\begin{align*}
\sup_{t \in I}\|\bm{\eta}(t)-\bm{y}(t)\|
  &= \sup_{t \in I}\|R(t)\bm{\phi}(t)-R(t)\bm{x}(t)\| \le \sup_{t \in I}\|R(t)\|\sup_{t \in I}\|\bm{\phi}(t)-\bm{x}(t)\|\\
  &\le \left(\sup_{t \in I}\|R(t)\|\right)\left(\sup_{t \in I}\|R^{-1}(t)\|\right)K\varepsilon.
\end{align*}
Hence, \eqref{Jordan} is Ulam stable on $I$, with an Ulam constant $\left(\sup_{t \in I}\|R(t)\|\right)\left(\sup_{t \in I}\|R^{-1}(t)\|\right)K$.
\end{proof}

Using the same method of the proof of Theorem \ref{uimpliesJ}, we also obtain the following result.

\begin{theorem}\label{Jimpliesu}
Suppose that $\|R(t)\|$ and $\|R^{-1}(t)\|$ are bounded on $I$. If equation \eqref{Jordan} is Ulam stable on $I$, with an Ulam constant $K$, then equation \eqref{unperturbed} is Ulam stable on $I$, and an Ulam constant is $\left(\sup_{t \in I}\|R(t)\|\right)\left(\sup_{t \in I}\|R^{-1}(t)\|\right)K$.
\end{theorem}

\begin{proof}
Assume that $\|R(t)\|$ and $\|R^{-1}(t)\|$ are bounded on $I$, and equation \eqref{Jordan} is Ulam stable on $I$, with an Ulam constant $K$. 
Let $\varepsilon >0$, and let $\bm{\phi}: I \to \C^n$ satisfy
$$ \sup_{t \in I}\|\bm{\phi}'-A(t)\bm{\phi}\| \le \varepsilon. $$
Let $\bm{\eta} = R(t)\bm{\phi}$. Then, we obtain
\begin{align*}
\sup_{t \in I}\|\bm{\eta}'-J(t)\bm{\eta}\|
  &= \sup_{t \in I}\|R(t)\bm{\phi}'+R'(t)\bm{\phi}-\left(R'(t)+R(t)A(t)\right)\bm{\phi}\| \\
  &= \sup_{t \in I}\|R(t)(\bm{\phi}'-A(t)\bm{\phi})\| \\
  &\le \varepsilon\sup_{t \in I}\|R(t)\|.
\end{align*}
Since $\|R(t)\|$ is bounded on $I$, and \eqref{Jordan} is Ulam stable on $I$, with an Ulam constant $K$, we see that there is a solution $\bm{y}: I \to \C^n$ of \eqref{Jordan} such that
\begin{equation}
 \sup_{t \in I}\|\bm{\eta}(t)-\bm{y}(t)\| \le K \varepsilon\sup_{t \in I}\|R(t)\| < \infty.
 \label{eq:2.03}
\end{equation}
Letting $\bm{x} = R^{-1}(t)\bm{y}$, we have
$$ \bm{x}' = R^{-1}(t)\bm{y}'+(R^{-1}(t))'\bm{y} = A(t)\bm{x}, $$
and so that $\bm{x}$ is a solution to \eqref{unperturbed}. In addition, by \eqref{eq:2.03} and the boundedness of $\|R^{-1}(t)\|$, we obtain
\begin{align*}
\sup_{t \in I}\|\bm{\phi}(t)-\bm{x}(t)\|
  &= \sup_{t \in I}\|R^{-1}(t)\bm{\eta}(t)-R^{-1}(t)\bm{y}(t)\| \le \sup_{t \in I}\|R^{-1}(t)\|\sup_{t \in I}\|\bm{\eta}(t)-\bm{y}(t)\|\\
  &\le \left(\sup_{t \in I}\|R(t)\|\right)\left(\sup_{t \in I}\|R^{-1}(t)\|\right)K\varepsilon.
\end{align*}
Hence, \eqref{unperturbed} is Ulam stable on $I$, with an Ulam constant $\left(\sup_{t \in I}\|R(t)\|\right)\left(\sup_{t \in I}\|R^{-1}(t)\|\right)K$.
\end{proof}

Combining Theorems \ref{iff1}--\ref{Jimpliesu}, we obtain the following result.

\begin{theorem}\label{iff2}
Suppose that $\|R(t)\|$ and $\|R^{-1}(t)\|$ are bounded on $I$. Then equation \eqref{perturbed} is Ulam stable on $I$ if and only if equation \eqref{Jordan} is Ulam stable on $I$.
\end{theorem}

\begin{proof}
We omit the proof because it is obvious.
\end{proof}


\section{Ulam stability of generalized Jordan normal form (I)}

We consider the two-dimensional linear differential system
\begin{equation}
 \bm{x}' = A(t)\bm{x}, \quad 
A(t) = 
\begin{pmatrix}
\lambda_1(t) & 0\\
0 & \lambda_2(t)
\end{pmatrix},
 \label{main1}
\end{equation}
where $\lambda_1$, $\lambda_2: I \to \C$ are continuous, and $\bm{x}\in \C^2$. In this case, we can easily find a fundamental matrix for this system. For example, a fundamental matrix of \eqref{main1} and its inverse matrix are as follows:
\begin{equation}
 X(t) = \begin{pmatrix}
e^{\int_{t_0}^t\lambda_1(s)ds} & 0\\
0 & e^{\int_{t_0}^t\lambda_2(s)ds}
\end{pmatrix}, \quad 
 X^{-1}(t) = \begin{pmatrix}
e^{-\int_{t_0}^t\lambda_1(s)ds} & 0\\
0 & e^{-\int_{t_0}^t\lambda_2(s)ds}
\end{pmatrix}
 \label{funda1}
\end{equation}
for $t \in I$, where $t_0 \in I$. In this section, we use the norm of the vector $\bm{v} = (v_1,v_2)^T$ as the maximum norm $\|\bm{v}\|_{\infty}:= \max\{|v_1|,|v_2|\}$. From the properties of the norm, if Ulam stability is shown for a norm, Ulam stability is guaranteed for any norm. Throughout this paper, the real part of $z \in \C$ will be written as $\Re(z)$. The first result in this section is as follows. 


\begin{theorem}\label{UlamS1}
Let $I$ be either $(a,b)$, $(a,b]$, $[a,b)$ or $[a,b]$, where $-\infty \le a < b \le \infty$. Then the following (i), (ii), and (iii) below hold:
\begin{itemize}
  \item[(i)] if
\begin{equation}
 \kappa_{11}(t):= \int_t^b e^{-\int_t^s\min\{\Re(\lambda_1(\tau)),\Re(\lambda_2(\tau))\}d\tau}ds \;\; \text{exists for all}\;\; t \in I,
 \label{kappa11}
\end{equation}
and $\sup_{t \in I}\kappa_{11}(t) < \infty$, then \eqref{main1} is Ulam stable on $I$, with an Ulam constant $K_{11} := \sup_{t \in I}\kappa_{11}(t)$. 
Furthermore, if
\begin{equation}
 \lim_{t\to b-0}\int_{t_0}^t\min\{\Re(\lambda_1(s)),\Re(\lambda_2(s))\}ds = \infty, \quad t_0 \in (a,b),
 \label{unbound11}
\end{equation}
then for any $\varepsilon>0$ and for any continuously differentiable function $\bm{\phi}(t)$ satisfying
\begin{equation}
 \sup_{t \in I}\|\bm{\phi}'(t) - A(t)\bm{\phi}(t)\|_{\infty} \le \varepsilon,
 \label{inequ}
\end{equation}
there exists the unique solution of \eqref{main1} denoted by
\[ \bm{x}(t) = X(t)\left(\lim_{t\to b-0} X^{-1}(t)\bm{\phi}(t)\right) \]
such that $\sup_{t \in I}\|\bm{\phi}(t)-\bm{x}(t)\|_{\infty} \le K_{11}\varepsilon$, where $X(t)$ and $X^{-1}(t)$ are given in \eqref{funda1}, and $\left(\lim_{t\to b-0} X^{-1}(t)\bm{\phi}(t)\right)$ is a well-defined constant vector;
  \item[(ii)] if
\begin{equation}
 \kappa_{12}(t):= \int_{a}^t e^{\int_s^t\max\{\Re(\lambda_1(\tau)),\Re(\lambda_2(\tau))\}d\tau}ds \;\; \text{exists for all}\;\; t \in I,
 \label{kappa12}
\end{equation}
and $\sup_{t \in I}\kappa_{12}(t) < \infty$, then \eqref{main1} is Ulam stable on $I$, with an Ulam constant $K_{12} := \sup_{t \in I}\kappa_{12}(t)$. 
Furthermore, if
\begin{equation}
 \lim_{t\to a+0}\int_t^{t_0}\min\{-\Re(\lambda_1(s)),-\Re(\lambda_2(s))\}ds = \infty, \quad t_0 \in (a,b),
 \label{unbound12}
\end{equation}
then for any $\varepsilon>0$ and for any continuously differentiable function $\bm{\phi}(t)$ satisfying \eqref{inequ}, there exists the unique solution of \eqref{main1} denoted by
\[ \bm{x}(t) = X(t)\left(\lim_{t\to a+0} X^{-1}(t)\bm{\phi}(t)\right) \]
such that $\sup_{t \in I}\|\bm{\phi}(t)-\bm{x}(t)\|_{\infty} \le K_{12}\varepsilon$, where $X(t)$ and $X^{-1}(t)$ are given in \eqref{funda1}, and $\left(\lim_{t\to a+0} X^{-1}(t)\bm{\phi}(t)\right)$ is a well-defined constant vector;
  \item[(iii)] if
\begin{equation}
 \kappa_{13}(t):= \max\left\{\int_t^b e^{-\int_t^s\Re(\lambda_1(\tau))d\tau}ds, \int_{a}^t e^{\int_s^t\Re(\lambda_2(\tau))d\tau}ds\right\} \;\; \text{exists for all}\;\; t \in I,
 \label{kappa13}
\end{equation}
and $\sup_{t \in I}\kappa_{13}(t) < \infty$, then \eqref{main1} is Ulam stable on $I$, with an Ulam constant $K_{13} := \sup_{t \in I}\kappa_{13}(t)$. 
Furthermore, if
\begin{equation}
 \lim_{t\to b-0}\int_{t_0}^t\Re(\lambda_1(s))ds = \lim_{t\to a+0}\int_{t_0}^t\Re(\lambda_2(s))ds = \infty, \quad t_0 \in (a,b),
 \label{unbound13}
\end{equation}
then for any $\varepsilon>0$ and for any continuously differentiable function $\bm{\phi}(t) = (\phi_1(t),\phi_2(t))^T$ satisfying \eqref{inequ}, there exists the unique solution of \eqref{main1} denoted by
\[ \bm{x}(t) = X(t)\begin{pmatrix}
\lim_{t\to b-0} \phi_1(t)e^{-\int_{t_0}^t\lambda_1(s)ds}\\
\lim_{t\to a+0} \phi_2(t)e^{\int_t^{t_0}\lambda_2(s)ds}
\end{pmatrix} \]
such that $\sup_{t \in I}\|\bm{\phi}(t)-\bm{x}(t)\|_{\infty} \le K_{13}\varepsilon$, where $\begin{pmatrix}
\lim_{t\to b-0} \phi_1(t)e^{-\int_{t_0}^t\lambda_1(s)ds}\\
\lim_{t\to a+0} \phi_2(t)e^{\int_t^{t_0}\lambda_2(s)ds}
\end{pmatrix}$ is a well-defined constant vector.
\end{itemize}
\end{theorem}

\begin{proof}
Case (i). Suppose that \eqref{kappa11} and $\sup_{t \in I}\kappa_{11}(t) < \infty$ hold. Let an arbitrary $\varepsilon>0$ be given, and let $\bm{\phi}(t)$ be continuously differentiable on $I$ and satisfy \eqref{inequ}. Define
\begin{equation}
 \bm{f}(t):= \bm{\phi}'(t) - A(t)\bm{\phi}(t), \quad t \in I.
 \label{deff}
\end{equation}
Then $\sup_{t \in I}\|\bm{f}(t)\|_{\infty} \le \varepsilon$ holds, and $\bm{\phi}(t)$ is a solution of \eqref{perturbed}. Thus, by the variation of parameters formula, we have
\begin{equation}
 \bm{\phi}(t):= X(t)X^{-1}(t_0)\bm{\phi}(t_0) + X(t)\int_{t_0}^tX^{-1}(s)\bm{f}(s)ds, \quad t_0 \in (a,b)
 \label{variation}
\end{equation}
for all $t \in I$, where $X(t)$ and $X^{-1}(t)$ are given in \eqref{funda1}. 
Now we consider a solution $\bm{x}(t)$ of \eqref{main1} defined by
\[ \bm{x}(t) := X(t)\bm{x}_{11}, \]
where
\begin{equation}
 \bm{x}_{11} := X^{-1}(t_0)\bm{\phi}(t_0)+\int_{t_0}^{b}X^{-1}(s)\bm{f}(s)ds.
 \label{x_11}
\end{equation}
Note that the (improper) integral in \eqref{x_11} converges. In fact, using \eqref{funda1}, \eqref{kappa11}, and the definition of the maximum norm, 
\begin{align*}
 \left\|\int_{t_0}^{b}X^{-1}(s)\bm{f}(s)ds\right\|_{\infty}
  &\le  \int_{t_0}^{b}\left\|X^{-1}(s)\bm{f}(s)\right\|_{\infty}ds\\
  &= \int_{t_0}^{b}\left\|
\begin{pmatrix}
e^{-\int_{t_0}^s\lambda_1(\tau)d\tau} & 0\\
0 & e^{-\int_{t_0}^s\lambda_2(\tau)d\tau}
\end{pmatrix}
\begin{pmatrix}
f_1(s)\\
f_2(s)
\end{pmatrix}
\right\|_{\infty}ds \\
  &= \int_{t_0}^{b}\left\|
\begin{pmatrix}
f_1(s)e^{-\int_{t_0}^s\lambda_1(\tau)d\tau}\\
f_2(s)e^{-\int_{t_0}^s\lambda_2(\tau)d\tau}
\end{pmatrix}
\right\|_{\infty}ds \\
  &= \int_{t_0}^{b}\max\left\{|f_1(s)|e^{-\int_{t_0}^s\Re(\lambda_1(\tau))d\tau}, |f_2(s)|e^{-\int_{t_0}^s\Re(\lambda_2(\tau))d\tau}\right\}ds \\
  &\le \int_{t_0}^{b} \|\bm{f}(s)\|_{\infty}e^{-\int_{t_0}^s\min\{\Re(\lambda_1(\tau)),\Re(\lambda_2(\tau))\}d\tau}ds \\
  &\le \varepsilon \int_{t_0}^{b}e^{-\int_{t_0}^s\min\{\Re(\lambda_1(\tau)),\Re(\lambda_2(\tau))\}d\tau}ds = \varepsilon \kappa_{11}(t_0) < \infty,
\end{align*}
where $f_1(t)$ and $f_2(t)$ are the components of the vector $\bm{f}(t)$; that is, $\bm{f}(t) = (f_1(t), f_2(t))^{T}$. Therefore, the constant vector $\bm{x}_{11}$ is well-defined, and $\bm{x}(t)$ is a solution of \eqref{main1}. By \eqref{variation} and \eqref{x_11}, we obtain
\begin{equation}
 \bm{\phi}(t)-\bm{x}(t) = -X(t)\int_t^{b}X^{-1}(s)\bm{f}(s)ds
 \label{phiminus11}
\end{equation}
for all $t \in I$. Hence, using \eqref{funda1}, \eqref{kappa11}, \eqref{phiminus11}, and the definition of the maximum norm again, we have
\begin{align*}
 \|\bm{\phi}(t)-\bm{x}(t)\|_{\infty}
  &\le  \int_t^{b}\left\|X(t)X^{-1}(s)\bm{f}(s)\right\|_{\infty}ds\\
  &= \int_t^{b}\left\|
\begin{pmatrix}
e^{-\int_t^s\lambda_1(\tau)d\tau} & 0\\
0 & e^{-\int_t^s\lambda_2(\tau)d\tau}
\end{pmatrix}
\begin{pmatrix}
f_1(s)\\
f_2(s)
\end{pmatrix}
\right\|_{\infty}ds \\
  &= \int_t^{b}\left\|
\begin{pmatrix}
f_1(s)e^{-\int_t^s\lambda_1(\tau)d\tau}\\
f_2(s)e^{-\int_t^s\lambda_2(\tau)d\tau}
\end{pmatrix}
\right\|_{\infty}ds \\
  &= \int_t^{b}\max\left\{|f_1(s)|e^{-\int_t^s\Re(\lambda_1(\tau))d\tau}, |f_2(s)|e^{-\int_t^s\Re(\lambda_2(\tau))d\tau}\right\}ds \\
  &\le \int_t^{b} \|\bm{f}(s)\|_{\infty}e^{-\int_t^s\min\{\Re(\lambda_1(\tau)),\Re(\lambda_2(\tau))\}d\tau}ds \\
  &\le \varepsilon \int_t^{b}e^{-\int_t^s\min\{\Re(\lambda_1(\tau)),\Re(\lambda_2(\tau))\}d\tau}ds = \varepsilon \kappa_{11}(t)
\end{align*}
for all $t \in I$, and thus,
\[ \sup_{t \in I}\|\bm{\phi}(t)-\bm{x}(t)\|_{\infty} \le \varepsilon \sup_{t \in I}\kappa_{11}(t) = K_{11}\varepsilon < \infty. \]
Therefore, \eqref{main1} is Ulam stable on $I$, with an Ulam constant $K_{11}$. Moreover, from \eqref{variation} and \eqref{x_11}, we have $\bm{x}_{11} = \lim_{t\to b-0} X^{-1}(t)\bm{\phi}(t)$, and so that $\bm{x}(t)$ is given by
\[ \bm{x}(t) = X(t)\left(\lim_{t\to b-0} X^{-1}(t)\bm{\phi}(t)\right). \]

Next, we discuss the uniqueness of the solution $\bm{x}(t)$. To show by contradiction, we assume that there is a solution $\bm{y}(t)$ of \eqref{main1} such that
\[ \bm{y}(t) := X(t)\bm{y}_{11}, \quad \bm{y}_{11} \ne \bm{x}_{11} \]
on $I$, and
\[ \sup_{t \in I}\|\bm{\phi}(t)-\bm{y}(t)\|_{\infty} \le K_{11}\varepsilon. \]
This implies
\[ \|X(t)(\bm{y}_{11}-\bm{x}_{11})\|_{\infty} \le \|\bm{\phi}(t)-\bm{y}(t)\|_{\infty} + \|\bm{\phi}(t)-\bm{x}(t)\|_{\infty} \le 2K_{11} \varepsilon \]
for all $t \in I$. However, by \eqref{unbound11}, we see that
\begin{align*}
 \lim_{t\to b-0}\|X(t)(\bm{y}_{11}-\bm{x}_{11})\|_{\infty}
  &= \lim_{t\to b-0}\left\|
\begin{pmatrix}
e^{\int_{t_0}^t\lambda_1(s)ds} & 0\\
0 & e^{\int_{t_0}^t\lambda_2(s)ds}
\end{pmatrix}
\begin{pmatrix}
y_1\\
y_2
\end{pmatrix}
\right\|_{\infty} \\
  &= \lim_{t\to b-0}\left\|
\begin{pmatrix}
y_1e^{\int_{t_0}^t\lambda_1(s)ds}\\
y_2e^{\int_{t_0}^t\lambda_2(s)ds}
\end{pmatrix}
\right\|_{\infty} \\
  &= \lim_{t\to b-0}\max\left\{|y_1|e^{\int_{t_0}^t\Re(\lambda_1(s))ds}, |y_2|e^{\int_{t_0}^t\Re(\lambda_2(s))ds}\right\} \\
  &\ge \max\{|y_1|, |y_2|\}\lim_{t\to b-0}e^{\int_{t_0}^t\min\{\Re(\lambda_1(s)),\Re(\lambda_2(s))\}ds} = \infty,
\end{align*}
where $\bm{y}_{11}-\bm{x}_{11} = (y_1,y_2)^T$ and $(y_1,y_2) \ne (0,0)$. This is a contradiction. Hence $\bm{x}(t)$ is the unique solution of \eqref{main1} satisfying $\sup_{t \in I}\|\bm{\phi}(t)-\bm{x}(t)\|_{\infty} \le K_{11}\varepsilon$.

Case (ii). Suppose that \eqref{kappa12} and $\sup_{t \in I}\kappa_{12}(t) < \infty$ hold. Let $\varepsilon>0$, and let $\bm{\phi}(t)$ be continuously differentiable on $I$ and satisfy \eqref{inequ}. Define $\bm{f}(t)$ by \eqref{deff}. Then $\sup_{t \in I}\|\bm{f}(t)\|_{\infty} \le \varepsilon$ holds, and $\bm{\phi}(t)$ satisfies \eqref{variation} for all $t \in I$, where $X(t)$ and $X^{-1}(t)$ are given in \eqref{funda1}. 
Now we consider a solution $\bm{x}(t)$ of \eqref{main1} defined by
\[ \bm{x}(t) := X(t)\bm{x}_{12}, \]
where
\begin{equation}
 \bm{x}_{12} := X^{-1}(t_0)\bm{\phi}(t_0)-\int_{a}^{t_0}X^{-1}(s)\bm{f}(s)ds.
 \label{x_12}
\end{equation}
We will prove that the (improper) integral in \eqref{x_12} converges. Put $\bm{f}(t) = (f_1(t), f_2(t))^{T}$. By \eqref{funda1}, \eqref{kappa12}, and the definition of the maximum norm, 
\begin{align*}
 \left\|\int_{a}^{t_0}X^{-1}(s)\bm{f}(s)ds\right\|_{\infty}
  &\le \int_{a}^{t_0} \left\|
\begin{pmatrix}
f_1(s)e^{\int_s^{t_0}\lambda_1(\tau)d\tau}\\
f_2(s)e^{\int_s^{t_0}\lambda_2(\tau)d\tau}
\end{pmatrix}
\right\|_{\infty}ds \\
  &\le \varepsilon \int_{a}^{t_0} e^{\int_s^{t_0}\max\{\Re(\lambda_1(\tau)),\Re(\lambda_2(\tau))\}d\tau}ds = \varepsilon \kappa_{12}(t_0) < \infty.
\end{align*}
Thus, $\bm{x}_{12}$ is well-defined. From \eqref{variation} and \eqref{x_12}, we obtain
\begin{equation}
 \bm{\phi}(t)-\bm{x}(t) = X(t)\int_{a}^t X^{-1}(s)\bm{f}(s)ds
 \label{phiminus12}
\end{equation}
for all $t \in I$. Hence, using \eqref{funda1}, \eqref{kappa12}, \eqref{phiminus12}, and the definition of the maximum norm again, we have
\begin{align*}
 \|\bm{\phi}(t)-\bm{x}(t)\|_{\infty}
  &\le \int_{a}^t \left\|
\begin{pmatrix}
f_1(s)e^{\int_s^t\lambda_1(\tau)d\tau}\\
f_2(s)e^{\int_s^t\lambda_2(\tau)d\tau}
\end{pmatrix}
\right\|_{\infty}ds \\
  &\le \varepsilon \int_{a}^t e^{\int_s^t\max\{\Re(\lambda_1(\tau)),\Re(\lambda_2(\tau))\}d\tau}ds = \varepsilon \kappa_{12}(t)
\end{align*}
for all $t \in I$, and thus,
\[ \sup_{t \in I}\|\bm{\phi}(t)-\bm{x}(t)\|_{\infty} \le \varepsilon \sup_{t \in I}\kappa_{12}(t) = K_{12}\varepsilon < \infty. \]
Consequently, \eqref{main1} is Ulam stable on $I$, with an Ulam constant $K_{12}$. 
Moreover, from \eqref{variation} and \eqref{x_12}, we have $\bm{x}_{12} = \lim_{t\to a+0} X^{-1}(t)\bm{\phi}(t)$, and so that $\bm{x}(t)$ is given by
\[ \bm{x}(t) = X(t)\left(\lim_{t\to a+0} X^{-1}(t)\bm{\phi}(t)\right). \]

Next, we will prove the uniqueness of the solution $\bm{x}(t)$. To show by contradiction, we assume that there is a solution $\bm{y}(t)$ of \eqref{main1} such that
\[ \bm{y}(t) := X(t)\bm{y}_{12}, \quad \bm{y}_{12} \ne \bm{x}_{12} \]
on $I$, and
\[ \sup_{t \in I}\|\bm{\phi}(t)-\bm{y}(t)\|_{\infty} \le K_{12}\varepsilon. \]
This implies
\[ \|X(t)(\bm{y}_{12}-\bm{x}_{12})\|_{\infty} \le \|\bm{\phi}(t)-\bm{y}(t)\|_{\infty} + \|\bm{\phi}(t)-\bm{x}(t)\|_{\infty} \le 2K_{12} \varepsilon \]
for all $t \in I$. However, by \eqref{unbound12}, we see that
\begin{align*}
 \lim_{t\to a+0}\|X(t)(\bm{y}_{12}-\bm{x}_{12})\|_{\infty}
  &= \lim_{t\to a+0}\max\left\{|y_1|e^{-\int_t^{t_0}\Re(\lambda_1(s))ds}, |y_2|e^{-\int_t^{t_0}\Re(\lambda_2(s))ds}\right\} \\
  &\ge \max\{|y_1|, |y_2|\}\lim_{t\to a+0}e^{\int_t^{t_0}\min\{-\Re(\lambda_1(s)),-\Re(\lambda_2(s))\}ds} = \infty,
\end{align*}
where $\bm{y}_{12}-\bm{x}_{12} = (y_1,y_2)^T$ and $(y_1,y_2) \ne (0,0)$. This is a contradiction. Hence $\bm{x}(t)$ is the unique solution of \eqref{main1} satisfying $\sup_{t \in I}\|\bm{\phi}(t)-\bm{x}(t)\|_{\infty} \le K_{12}\varepsilon$. 

Case (iii). Suppose that \eqref{kappa13} and $\sup_{t \in I}\kappa_{13}(t) < \infty$ hold. Let $\varepsilon>0$, and let $\bm{\phi}(t)$ be continuously differentiable on $I$ and satisfy \eqref{inequ}. Define $\bm{f}(t)$ by \eqref{deff}. Then $\sup_{t \in I}\|\bm{f}(t)\|_{\infty} \le \varepsilon$ holds, and $\bm{\phi}(t)$ satisfies \eqref{variation} for all $t \in I$, where $X(t)$ and $X^{-1}(t)$ are given in \eqref{funda1}. 
Put $\bm{f}(t) = (f_1(t), f_2(t))^{T}$. 
Now we consider a solution $\bm{x}(t)$ of \eqref{main1} defined by
\[ \bm{x}(t) := X(t)\bm{x}_{13}, \]
where
\begin{equation}
 \bm{x}_{13} := X^{-1}(t_0)\bm{\phi}(t_0) + 
\begin{pmatrix}
\int_{t_0}^{b} f_1(s)e^{-\int_{t_0}^s\lambda_1(\tau)d\tau}ds\\
- \int_{a}^{t_0} f_2(s)e^{\int_s^{t_0}\lambda_2(\tau)d\tau}ds
\end{pmatrix}
.
 \label{x_13}
\end{equation}
We will prove that the (improper) integrals in \eqref{x_13} converge. By \eqref{funda1}, \eqref{kappa13}, and the definition of the maximum norm, 
\begin{align*}
 \left\|
\begin{pmatrix}
\int_{t_0}^{b} f_1(s)e^{-\int_{t_0}^s\lambda_1(\tau)d\tau}ds\\
- \int_{a}^{t_0} f_2(s)e^{\int_s^{t_0}\lambda_2(\tau)d\tau}ds
\end{pmatrix}
\right\|_{\infty} 
  &= \max\left\{\left|\int_{t_0}^{b} f_1(s)e^{-\int_{t_0}^s\lambda_1(\tau)d\tau}ds\right|, \left|\int_{a}^{t_0} f_2(s)e^{\int_s^{t_0}\lambda_2(\tau)d\tau}ds\right|\right\} \\
  &\le \max\left\{\int_{t_0}^{b} |f_1(s)|e^{-\int_{t_0}^s\Re(\lambda_1(\tau))d\tau}ds, \int_{a}^{t_0} |f_2(s)|e^{\int_s^{t_0}\Re(\lambda_2(\tau))d\tau}ds\right\} \\
  &\le \max\left\{\int_{t_0}^{b} \|\bm{f}(s)\|_{\infty}e^{-\int_{t_0}^s\Re(\lambda_1(\tau))d\tau}ds, \int_{a}^{t_0} \|\bm{f}(s)\|_{\infty}e^{\int_s^{t_0}\Re(\lambda_2(\tau))d\tau}ds\right\} \\
  &\le \varepsilon \max\left\{\int_{t_0}^{b} e^{-\int_{t_0}^s\Re(\lambda_1(\tau))d\tau}ds, \int_{a}^{t_0} e^{\int_s^{t_0}\Re(\lambda_2(\tau))d\tau}ds \right\} \\
  &= \varepsilon \kappa_{13}(t_0) < \infty.
\end{align*}
Thus, $\bm{x}_{13}$ is well-defined. From \eqref{variation} and \eqref{x_13}, we obtain
\begin{align}
 \bm{\phi}(t)-\bm{x}(t) 
  &= X(t)\left(\int_{t_0}^tX^{-1}(s)\bm{f}(s)ds - \begin{pmatrix}
\int_{t_0}^{b} f_1(s)e^{-\int_{t_0}^s\lambda_1(\tau)d\tau}ds\\
- \int_{a}^{t_0} f_2(s)e^{\int_s^{t_0}\lambda_2(\tau)d\tau}ds
\end{pmatrix}\right) \nonumber\\
 &= 
\begin{pmatrix}
-\int_t^{b} f_1(s)e^{-\int_t^s\lambda_1(\tau)d\tau}ds\\
\int_{a}^t f_2(s)e^{\int_s^t\lambda_2(\tau)d\tau}ds
\end{pmatrix}
 \label{phiminus13}
\end{align}
for all $t \in I$. Hence, using \eqref{funda1}, \eqref{kappa13}, \eqref{phiminus13}, and the definition of the maximum norm again, we have
\begin{align*}
 \|\bm{\phi}(t)-\bm{x}(t)\|_{\infty}
  &= \max\left\{\left|\int_t^{b} f_1(s)e^{-\int_t^s\lambda_1(\tau)d\tau}ds\right|, \left|\int_{a}^t f_2(s)e^{\int_s^t\lambda_2(\tau)d\tau}ds\right|\right\} \\
  &\le \max\left\{\int_t^{b} |f_1(s)|e^{-\int_t^s\Re(\lambda_1(\tau))d\tau}ds, \int_{a}^t |f_2(s)|e^{\int_s^t\Re(\lambda_2(\tau))d\tau}ds\right\} \\
  &\le \max\left\{\int_t^{b} \|\bm{f}(s)\|_{\infty}e^{-\int_t^s\Re(\lambda_1(\tau))d\tau}ds, \int_{a}^t \|\bm{f}(s)\|_{\infty}e^{\int_s^t\Re(\lambda_2(\tau))d\tau}ds\right\} \\
  &\le \varepsilon \max\left\{\int_t^{b} e^{-\int_t^s\Re(\lambda_1(\tau))d\tau}ds, \int_{a}^t e^{\int_s^t\Re(\lambda_2(\tau))d\tau}ds \right\} \\
  &= \varepsilon \kappa_{13}(t)
\end{align*}
for all $t \in I$, and thus,
\[ \sup_{t \in I}\|\bm{\phi}(t)-\bm{x}(t)\|_{\infty} \le \varepsilon \sup_{t \in I}\kappa_{13}(t) = K_{13}\varepsilon < \infty. \]
Consequently, \eqref{main1} is Ulam stable on $I$, with an Ulam constant $K_{13}$. 
Moreover, from \eqref{variation} and \eqref{x_13}, we have
\[ \bm{x}_{13} = \begin{pmatrix}
\lim_{t\to b-0} \phi_1(t)e^{-\int_{t_0}^t\lambda_1(s)ds}\\
\lim_{t\to a+0} \phi_2(t)e^{\int_t^{t_0}\lambda_2(s)ds}
\end{pmatrix}, \]
where $\bm{\phi}(t) = (\phi_1(t),\phi_2(t))^T$. Hence $\bm{x}(t)$ is given by
\[ \bm{x}(t) = X(t)\begin{pmatrix}
\lim_{t\to b-0} \phi_1(t)e^{-\int_{t_0}^t\lambda_1(s)ds}\\
\lim_{t\to a+0} \phi_2(t)e^{\int_t^{t_0}\lambda_2(s)ds}
\end{pmatrix}. \]

Next, we will prove the uniqueness of the solution $\bm{x}(t)$. To show by contradiction, we assume that there is a solution $\bm{y}(t)$ of \eqref{main1} such that
\[ \bm{y}(t) := X(t)\bm{y}_{13}, \quad \bm{y}_{13} \ne \bm{x}_{13} \]
on $I$, and
\[ \sup_{t \in I}\|\bm{\phi}(t)-\bm{y}(t)\|_{\infty} \le K_{13}\varepsilon. \]
This implies
\[ \|X(t)(\bm{y}_{13}-\bm{x}_{13})\|_{\infty} \le \|\bm{\phi}(t)-\bm{y}(t)\|_{\infty} + \|\bm{\phi}(t)-\bm{x}(t)\|_{\infty} \le 2K_{13} \varepsilon \]
for all $t \in I$. However, by \eqref{unbound13}, we see that
\begin{align*}
 \lim_{t\to b-0}\|X(t)(\bm{y}_{13}-\bm{x}_{13})\|_{\infty}
  &= \lim_{t\to b-0}\max\left\{|y_1|e^{\int_{t_0}^t\Re(\lambda_1(s))ds}, |y_2|e^{\int_{t_0}^t\Re(\lambda_2(s))ds}\right\} \\
  &\ge \lim_{t\to b-0}|y_1|e^{\int_{t_0}^t\Re(\lambda_1(s))ds} = \infty,
\end{align*}
if $y_1\ne 0$, and
\begin{align*}
 \lim_{t\to a+0}\|X(t)(\bm{y}_{13}-\bm{x}_{13})\|_{\infty}
  &= \lim_{t\to a+0}\max\left\{|y_1|e^{-\int_t^{t_0}\Re(\lambda_1(s))ds}, |y_2|e^{-\int_t^{t_0}\Re(\lambda_2(s))ds}\right\} \\
  &\ge \lim_{t\to a+0}|y_2|e^{-\int_t^{t_0}\Re(\lambda_2(s))ds} = \infty,
\end{align*}
if $y_2\ne 0$, where $\bm{y}_{13}-\bm{x}_{13} = (y_1,y_2)^T$ and $(y_1,y_2) \ne (0,0)$. This is a contradiction. Hence $\bm{x}(t)$ is the unique solution of \eqref{main1} satisfying $\sup_{t \in I}\|\bm{\phi}(t)-\bm{x}(t)\|_{\infty} \le K_{13}\varepsilon$. Thus, the proof is now complete.
\end{proof}

Next, we discuss the lower bound of the Ulam constants.


\begin{theorem}\label{lower1}
Let $I$ be either $(a,b)$, $(a,b]$, $[a,b)$ or $[a,b]$, where $a \le a < b \le \infty$. Then the following (i), (ii) and (iii) below hold:
\begin{itemize}
  \item[(i)] if \eqref{kappa11}, \eqref{unbound11}, and $\sup_{t \in I}\kappa_{11}(t) < \infty$ hold, then \eqref{main1} is Ulam stable on $I$, and any Ulam constant is greater than or equal to $K_{11} = \sup_{t \in I}\kappa_{11}(t)$;
  \item[(ii)] if \eqref{kappa12}, \eqref{unbound12}, and $\sup_{t \in I}\kappa_{12}(t) < \infty$ hold, then \eqref{main1} is Ulam stable on $I$, and any Ulam constant is greater than or equal to $K_{12} = \sup_{t \in I}\kappa_{12}(t)$;
  \item[(iii)] if \eqref{kappa13}, \eqref{unbound13}, and $\sup_{t \in I}\kappa_{13}(t) < \infty$ hold, then \eqref{main1} is Ulam stable on $I$, and any Ulam constant is greater than or equal to $K_{13} = \sup_{t \in I}\kappa_{13}(t)$.
\end{itemize}
\end{theorem}

\begin{proof}
Case (i). Suppose that \eqref{kappa11}, \eqref{unbound11}, and $\sup_{t \in I}\kappa_{11}(t) < \infty$ hold. Define
\begin{equation}
 \bm{f}_1(t) := \varepsilon \begin{pmatrix}
e^{i\int_{t_0}^t\Im(\lambda_1(s))ds}\\
e^{i\int_{t_0}^t\Im(\lambda_2(s))ds}
\end{pmatrix},\quad t \in I,
 \label{f1}
\end{equation}
where $\Im(z)$ is the imaginary part of $z \in \C$, and $i$ is the imaginary unit. Now we consider the solution $\bm{\phi}_1(t)$ of the linear differential system
\[ \bm{\phi}_1' - A(t)\bm{\phi}_1 = \bm{f}_1(t). \]
Since $\|\bm{f}_1(t)\|_{\infty} = \varepsilon$ holds for all $t \in I$, $\bm{\phi}_1(t)$ satisfies \eqref{inequ}. By Theorem~\ref{UlamS1}~(i), there exists the unique solution of \eqref{main1} denoted by
\[ \bm{x}_1(t) = X(t)\left(\lim_{t\to b-0} X^{-1}(t)\bm{\phi}_1(t)\right) \]
such that $\sup_{t \in I}\|\bm{\phi}_1(t)-\bm{x}_1(t)\|_{\infty} \le K_{11}\varepsilon$, where
\[ K_{11} = \sup_{t \in I}\kappa_{11}(t) = \sup_{t \in I}\int_t^{b}e^{-\int_t^s\min\{\Re(\lambda_1(\tau)),\Re(\lambda_2(\tau))\}d\tau}ds. \]
Hence, we obtain \eqref{phiminus11} (see the proof of Theorem~\ref{UlamS1} (i)), and so that
\begin{align*}
 \left\|\bm{\phi}_1(t)-\bm{x}_1(t)\right\|_{\infty} 
  &= \left\|X(t)\int_t^{b}X^{-1}(s)\bm{f}_1(s)ds\right\|_{\infty}
   = \varepsilon \left\|\begin{pmatrix}
\int_t^{b} e^{-\int_t^s\lambda_1(\tau)d\tau+i\int_{t_0}^s\Im(\lambda_1(\tau))d\tau}ds\\
\int_t^{b} e^{-\int_t^s\lambda_2(\tau)d\tau+i\int_{t_0}^s\Im(\lambda_2(\tau))d\tau}ds
\end{pmatrix}\right\|_{\infty}\\
  &= \varepsilon \left\|\begin{pmatrix}
 e^{-\int_t^{t_0}\lambda_1(\tau)d\tau}\int_t^{b} e^{-\int_{t_0}^s\Re(\lambda_1(\tau))d\tau}ds\\
 e^{-\int_t^{t_0}\lambda_2(\tau)d\tau}\int_t^{b} e^{-\int_{t_0}^s\Re(\lambda_2(\tau))d\tau}ds
\end{pmatrix}\right\|_{\infty}\\
  &= \varepsilon \max\left\{\left|e^{-\int_t^{t_0}\lambda_1(\tau)d\tau}\int_t^{b} e^{-\int_{t_0}^s\Re(\lambda_1(\tau))d\tau}ds\right|, 
 \left|e^{-\int_t^{t_0}\lambda_2(\tau)d\tau}\int_t^{b} e^{-\int_{t_0}^s\Re(\lambda_2(\tau))d\tau}ds\right|\right\}\\
  &= \varepsilon \max\left\{\left|e^{-\int_t^{t_0}\lambda_1(\tau)d\tau}\right|\int_t^{b} e^{-\int_{t_0}^s\Re(\lambda_1(\tau))d\tau}ds, 
 \left|e^{-\int_t^{t_0}\lambda_2(\tau)d\tau}\right|\int_t^{b} e^{-\int_{t_0}^s\Re(\lambda_2(\tau))d\tau}ds\right\}\\
  &= \varepsilon \max\left\{e^{-\int_t^{t_0}\Re(\lambda_1(\tau))d\tau} \int_t^{b} e^{-\int_{t_0}^s\Re(\lambda_1(\tau))d\tau}ds, 
 e^{-\int_t^{t_0}\Re(\lambda_2(\tau))d\tau} \int_t^{b} e^{-\int_{t_0}^s\Re(\lambda_2(\tau))d\tau}ds\right\}\\
  &= \varepsilon \max\left\{\int_t^{b} e^{-\int_t^s\Re(\lambda_1(\tau))d\tau}ds, 
 \int_t^{b} e^{-\int_t^s\Re(\lambda_2(\tau))d\tau}ds\right\}\\
  &= \varepsilon \kappa_{11}(t).
\end{align*}
This implies that
\[ \sup_{t \in I}\left\|\bm{\phi}_1(t)-\bm{x}_1(t)\right\|_{\infty} = K_{11}\varepsilon < \infty. \]
Since $\bm{x}_1(t)$ is the unique solution of \eqref{main1} satisfying $\sup_{t \in I}\|\bm{\phi}_1(t)-\bm{x}_1(t)\|_{\infty} \le K_{11}\varepsilon$, there is no solution $\bm{x}(t)$ to \eqref{main1} that satisfies $\sup_{t \in I}\|\bm{\phi}_1(t)-\bm{x}(t)\|_{\infty} < K_{11}\varepsilon$. Consequently, the Ulam constant is at least $K_{11}$. 

Case (ii). Suppose that \eqref{kappa12}, \eqref{unbound12}, and $\sup_{t \in I}\kappa_{12}(t) < \infty$ hold. Consider the solution $\bm{\phi}_2(t)$ of the linear differential system
\[ \bm{\phi}_2' - A(t)\bm{\phi}_2 = \bm{f}_1(t), \]
where $\bm{f}_1(t)$ is given by \eqref{f1}. 
Since $\|\bm{f}_1(t)\|_{\infty} = \varepsilon$ holds for all $t \in I$, $\bm{\phi}_2(t)$ satisfies \eqref{inequ}. By Theorem~\ref{UlamS1}~(ii), there exists the unique solution of \eqref{main1} denoted by
\[ \bm{x}_2(t) = X(t)\left(\lim_{t\to a+0} X^{-1}(t)\bm{\phi}_2(t)\right) \]
such that $\sup_{t \in I}\|\bm{\phi}_2(t)-\bm{x}_2(t)\|_{\infty} \le K_{12}\varepsilon$, where $K_{12} = \sup_{t \in I}\kappa_{12}(t)$. Hence, we obtain \eqref{phiminus12} (see the proof of Theorem~\ref{UlamS1} (ii)), and so that
\begin{align*}
 \left\|\bm{\phi}_2(t)-\bm{x}_2(t)\right\|_{\infty} &= \left\|X(t)\int_{a}^t X^{-1}(s)\bm{f}_1(s)ds\right\|_{\infty}
  = \varepsilon \left\|\begin{pmatrix}
 e^{\int_{t_0}^t\lambda_1(\tau)d\tau}\int_{a}^t e^{\int_s^{t_0}\Re(\lambda_1(\tau))d\tau}ds\\
 e^{\int_{t_0}^t\lambda_2(\tau)d\tau}\int_{a}^t e^{\int_s^{t_0}\Re(\lambda_2(\tau))d\tau}ds
\end{pmatrix}\right\|_{\infty}\\
  &= \varepsilon \max\left\{\left|e^{\int_{t_0}^t\lambda_1(\tau)d\tau}\right|\int_{a}^t e^{\int_s^{t_0}\Re(\lambda_1(\tau))d\tau}ds, 
 \left|e^{\int_{t_0}^t\lambda_2(\tau)d\tau}\right|\int_{a}^t e^{\int_s^{t_0}\Re(\lambda_2(\tau))d\tau}ds\right\}\\
  &= \varepsilon \max\left\{e^{\int_{t_0}^t\Re(\lambda_1(\tau))d\tau} \int_{a}^t e^{\int_s^{t_0}\Re(\lambda_1(\tau))d\tau}ds, 
 e^{\int_{t_0}^t\Re(\lambda_2(\tau))d\tau} \int_{a}^t e^{\int_s^{t_0}\Re(\lambda_2(\tau))d\tau}ds\right\}\\
  &= \varepsilon \kappa_{12}(t).
\end{align*}
Thus, we obtain $\sup_{t \in I}\left\|\bm{\phi}_2(t)-\bm{x}_2(t)\right\|_{\infty} = K_{12}\varepsilon < \infty$. Consequently, any Ulam constant is at least~$K_{12}$. 

Case (iii). Suppose that \eqref{kappa13}, \eqref{unbound13}, and $\sup_{t \in I}\kappa_{13}(t) < \infty$ hold. Consider the solution $\bm{\phi}_3(t)$ of the linear differential system
\[ \bm{\phi}_3' - A(t)\bm{\phi}_3 = \bm{f}_1(t), \]
where $\bm{f}_1(t)$ is given by \eqref{f1}. 
Since $\|\bm{f}_1(t)\|_{\infty} = \varepsilon$ holds for all $t \in I$, $\bm{\phi}_3(t)$ satisfies \eqref{inequ}. By Theorem~\ref{UlamS1}~(iii), there exists the unique solution of \eqref{main1} denoted by
\[ \bm{x}_3(t) = X(t)\begin{pmatrix}
\lim_{t\to b-0} \phi_1(t)e^{-\int_{t_0}^t\lambda_1(s)ds}\\
\lim_{t\to a+0} \phi_2(t)e^{\int_t^{t_0}\lambda_2(s)ds}
\end{pmatrix} \]
such that $\sup_{t \in I}\|\bm{\phi}_3(t)-\bm{x}_3(t)\|_{\infty} \le K_{13}\varepsilon$, where $K_{13} = \sup_{t \in I}\kappa_{13}(t)$. Hence, we obtain \eqref{phiminus13} (see the proof of Theorem~\ref{UlamS1} (iii)), and so that
\begin{align*}
 \left\|\bm{\phi}_3(t)-\bm{x}_3(t)\right\|_{\infty} &= \left\|\begin{pmatrix}
-\int_t^{b} f_1(s)e^{-\int_t^s\lambda_1(\tau)d\tau}ds\\
\int_{a}^t f_2(s)e^{\int_s^t\lambda_2(\tau)d\tau}ds
\end{pmatrix}\right\|_{\infty}
   = \varepsilon \left\|\begin{pmatrix}
-\int_t^{b} e^{-\int_t^s\lambda_1(\tau)d\tau+i\int_{t_0}^s\Im(\lambda_1(\tau))d\tau}ds\\
\int_{a}^t e^{\int_s^t\lambda_2(\tau)d\tau+i\int_{t_0}^s\Im(\lambda_2(\tau))d\tau}ds
\end{pmatrix}\right\|_{\infty}\\
  &= \varepsilon \left\|\begin{pmatrix}
-e^{-\int_t^{t_0}\lambda_1(\tau)d\tau}\int_t^{b} e^{-\int_{t_0}^s\Re(\lambda_1(\tau))d\tau}ds\\
 e^{\int_{t_0}^t\lambda_2(\tau)d\tau}\int_{a}^t e^{\int_s^{t_0}\Re(\lambda_2(\tau))d\tau}ds
\end{pmatrix}\right\|_{\infty}\\
  &= \varepsilon \max\left\{\left|e^{-\int_t^{t_0}\lambda_1(\tau)d\tau}\right|\int_t^{b} e^{-\int_{t_0}^s\Re(\lambda_1(\tau))d\tau}ds, 
 \left|e^{\int_{t_0}^t\lambda_2(\tau)d\tau}\right|\int_{a}^t e^{\int_s^{t_0}\Re(\lambda_2(\tau))d\tau}ds\right\}\\
  &= \varepsilon \max\left\{e^{-\int_t^{t_0}\Re(\lambda_1(\tau))d\tau} \int_t^{b} e^{-\int_{t_0}^s\Re(\lambda_1(\tau))d\tau}ds, 
 e^{\int_{t_0}^t\Re(\lambda_2(\tau))d\tau} \int_{a}^t e^{\int_s^{t_0}\Re(\lambda_2(\tau))d\tau}ds\right\}\\
  &= \varepsilon \kappa_{13}(t).
\end{align*}
Thus, we have $\sup_{t \in I}\left\|\bm{\phi}_3(t)-\bm{x}_3(t)\right\|_{\infty} = K_{13}\varepsilon < \infty$. Consequently, any Ulam constant is at least~$K_{13}$. The proof is now complete.
\end{proof}

Using Theorems \ref{UlamS1} and \ref{lower1}, we obtain the following result, immediately. 


\begin{theorem}\label{best1}
Let $I$ be either $(a,b)$, $(a,b]$, $[a,b)$ or $[a,b]$, where $a \le a < b \le \infty$. Then the following (i), (ii) and (iii) below hold:
\begin{itemize}
  \item[(i)] if \eqref{kappa11}, \eqref{unbound11}, and $\sup_{t \in I}\kappa_{11}(t) < \infty$ hold, then \eqref{main1} is Ulam stable on $I$, and the best Ulam constant is $K_{11} = \sup_{t \in I}\kappa_{11}(t)$;
  \item[(ii)] if \eqref{kappa12}, \eqref{unbound12}, and $\sup_{t \in I}\kappa_{12}(t) < \infty$ hold, then \eqref{main1} is Ulam stable on $I$, and the best Ulam constant is $K_{12} = \sup_{t \in I}\kappa_{12}(t)$;
  \item[(iii)] if \eqref{kappa13}, \eqref{unbound13}, and $\sup_{t \in I}\kappa_{13}(t) < \infty$ hold, then \eqref{main1} is Ulam stable on $I$, and the best Ulam constant is $K_{13} = \sup_{t \in I}\kappa_{13}(t)$.
\end{itemize}
\end{theorem}

When $\lambda_1(t) \equiv \lambda_1$ and $\lambda_2(t) \equiv \lambda_2$, \eqref{main1} reduces to the constant coefficients two-dimensional linear differential system
\begin{equation}
 \bm{x}' = A\bm{x}, \quad 
A = 
\begin{pmatrix}
\lambda_1 & 0\\
0 & \lambda_2
\end{pmatrix}.
 \label{const1}
\end{equation}
For this system, we have the following result.


\begin{corollary}\label{constbest1}
Let $I=\R$ and $\Re(\lambda_1)\Re(\lambda_2) \ne 0$. Then \eqref{const1} is Ulam stable on $\R$, and the best Ulam constant is $K_{c1} := \max\left\{\frac{1}{|\Re(\lambda_1)|}, \frac{1}{|\Re(\lambda_2)|}\right\}$. 
\end{corollary}

\begin{proof}
First, we consider the case $\Re(\lambda_1) \ge \Re(\lambda_2) > 0$. Since $\min\{\Re(\lambda_1),\Re(\lambda_2)\} = \Re(\lambda_2)$ holds, we have
\[ \kappa_{11}(t) = \int_t^{\infty} e^{-\int_t^s\min\{\Re(\lambda_1),\Re(\lambda_2)\}d\tau}ds = \int_t^{\infty} e^{-\Re(\lambda_2)(s-t) d\tau}ds = \frac{1}{\Re(\lambda_2)} = \max\left\{\frac{1}{|\Re(\lambda_1)|}, \frac{1}{|\Re(\lambda_2)|}\right\} \]
for all $t \in \R$. Using Theorem~\ref{best1} (i) with $I=\R$, \eqref{const1} is Ulam stable on $\R$, and the best Ulam constant is $K_{c1}$. 

Next, we consider the case $0 > \Re(\lambda_1) \ge \Re(\lambda_2)$. Since $\max\{\Re(\lambda_1),\Re(\lambda_2)\} = \Re(\lambda_1)$ holds, we have
\[ \kappa_{12}(t) = \int_{-\infty}^t e^{\int_s^t\max\{\Re(\lambda_1),\Re(\lambda_2)\}d\tau}ds = \int_{-\infty}^t e^{\Re(\lambda_1)(t-s)}ds = \frac{1}{-\Re(\lambda_1)} = \max\left\{\frac{1}{|\Re(\lambda_1)|}, \frac{1}{|\Re(\lambda_2)|}\right\} \]
for all $t \in \R$. Using Theorem~\ref{best1} (ii) with $I=\R$, \eqref{const1} is Ulam stable on $\R$, and the best Ulam constant is $K_{c1}$. 

Finally, we consider the case $\Re(\lambda_1) > 0 > \Re(\lambda_2)$. Then
\begin{align*}
 \kappa_{13}(t) &= \max\left\{\int_t^{\infty} e^{-\int_t^s\Re(\lambda_1)d\tau}ds, \int_{-\infty}^t e^{\int_s^t\Re(\lambda_2)d\tau}ds\right\}\\
  &= \max\left\{\frac{1}{\Re(\lambda_1)}, \frac{1}{-\Re(\lambda_2)}\right\} 
   = \max\left\{\frac{1}{|\Re(\lambda_1)|}, \frac{1}{|\Re(\lambda_2)|}\right\}
\end{align*}
for all $t \in \R$. By Theorem~\ref{best1} (iii) with $I=\R$, \eqref{const1} is Ulam stable on $\R$, and the best Ulam constant is $K_{c1}$. 
\end{proof}

\begin{remark}
When $\lambda_1$, $\lambda_2 \in \R$ and $\lambda_1\lambda_2>0$, we can check that the Ulam constant
\[ \max\left\{\frac{1}{|\Re(\lambda_1)|}, \frac{1}{|\Re(\lambda_2)|}\right\} = \max\left\{\frac{1}{|\lambda_1|}, \frac{1}{|\lambda_2|}\right\} = \|A^{-1}\|_{\infty} \]
in Corollary \ref{constbest1} is the best Ulam constant, by using the result in \cite[Theorem 4.3.]{AndOni5}, where $\|\cdot\|_{\infty}$ is the induced matrix norm defined by $\|M\|_{\infty} := \max_{j\in \{1,2\}}\{|m_{j1}|+|m_{j2}|\}$ for $M = \begin{pmatrix}
m_{11} & m_{12}\\
m_{21} & m_{22}
\end{pmatrix}$. However, to our knowledge no best Ulam constant for \eqref{const1} has ever been derived for the other cases satisfying $\Re(\lambda_1)\Re(\lambda_2) \ne 0$. In other words, it can be argued that new results have been obtained from this study, even if only for constant coefficients.
\end{remark}


\section{Ulam stability of generalized Jordan normal form (II)}

In this section, we consider the two-dimensional linear differential system
\begin{equation}
 \bm{x}' = A(t)\bm{x}, \quad 
A(t) = 
\begin{pmatrix}
\lambda(t) & \mu(t)\\
0 & \lambda(t)
\end{pmatrix},
 \label{main2}
\end{equation}
where $\lambda$, $\mu: I \to \C$ are continuous, and $\bm{x}\in \C^2$. A fundamental matrix of \eqref{main2} and its inverse matrix are as follows:
\begin{equation}
 X(t) = e^{\int_{t_0}^t\lambda(s)ds}\begin{pmatrix}
1 & \int_{t_0}^t\mu(s)ds\\
0 & 1
\end{pmatrix}, \quad 
 X^{-1}(t) = e^{-\int_{t_0}^t\lambda(s)ds}\begin{pmatrix}
1 & -\int_{t_0}^t\mu(s)ds\\
0 & 1
\end{pmatrix}
 \label{funda2}
\end{equation}
for $t \in I$, where $t_0 \in I$. In this section, the maximum norm $\|\cdot\|_{\infty}$ will be used. The obtained result is as follows. 


\begin{theorem}\label{UlamS2}
Let $I$ be either $(a,b)$, $(a,b]$, $[a,b)$ or $[a,b]$, where $-\infty \le a < b \le \infty$. Then the following (i) and (ii) below hold:
\begin{itemize}
  \item[(i)] if
\begin{equation}
 \kappa_{21}(t):= \int_t^{b} \left(1+\left|\int_t^s\mu(\tau)d\tau\right|\right)e^{-\int_t^s\Re(\lambda(\tau))d\tau}ds \;\; \text{exists for all}\;\; t \in I,
 \label{kappa21}
\end{equation}
and $\sup_{t \in I}\kappa_{21}(t) < \infty$, then \eqref{main2} is Ulam stable on $I$, with an Ulam constant $K_{21} := \sup_{t \in I}\kappa_{21}(t)$. 
Furthermore, if
\begin{equation}
 \lim_{t\to b-0}\int_{t_0}^t\Re(\lambda(s))ds = \infty, \quad t_0 \in (a,b),
 \label{unbound21}
\end{equation}
then for any $\varepsilon>0$ and for any continuously differentiable function $\bm{\phi}(t)$ satisfying \eqref{inequ}, there exists the unique solution of \eqref{main2} denoted by
\[ \bm{x}(t) = X(t)\left(\lim_{t\to b-0} X^{-1}(t)\bm{\phi}(t)\right) \]
such that $\sup_{t \in I}\|\bm{\phi}(t)-\bm{x}(t)\|_{\infty} \le K_{21}\varepsilon$, where $X(t)$ and $X^{-1}(t)$ are given in \eqref{funda2}, and $\left(\lim_{t\to b-0} X^{-1}(t)\bm{\phi}(t)\right)$ is a well-defined constant vector;
  \item[(ii)] if
\begin{equation}
 \kappa_{22}(t):= \int_{a}^t \left(1+\left|\int_s^t\mu(\tau)d\tau\right|\right)e^{\int_s^t\Re(\lambda(\tau))d\tau}ds \;\; \text{exists for all}\;\; t \in I,
 \label{kappa22}
\end{equation}
and $\sup_{t \in I}\kappa_{22}(t) < \infty$, then \eqref{main2} is Ulam stable on $I$, with an Ulam constant $K_{22} := \sup_{t \in I}\kappa_{22}(t)$. 
Furthermore, if
\begin{equation}
 \lim_{t\to a+0}\int_t^{t_0}\Re(\lambda(s))ds = -\infty, \quad t_0 \in (a,b),
 \label{unbound22}
\end{equation}
then for any $\varepsilon>0$ and for any continuously differentiable function $\bm{\phi}(t)$ satisfying \eqref{inequ}, there exists the unique solution of \eqref{main2} denoted by
\[ \bm{x}(t) = X(t)\left(\lim_{t\to a+0} X^{-1}(t)\bm{\phi}(t)\right) \]
such that $\sup_{t \in I}\|\bm{\phi}(t)-\bm{x}(t)\|_{\infty} \le K_{22}\varepsilon$, where $X(t)$ and $X^{-1}(t)$ are given in \eqref{funda2}, and $\left(\lim_{t\to a+0} X^{-1}(t)\bm{\phi}(t)\right)$ is a well-defined constant vector.
\end{itemize}
\end{theorem}

\begin{proof}
Since the proof of this theorem can be proved in the same way as Theorem \ref{UlamS1}, only the essential parts will be written below. 

Case (i). Suppose that \eqref{kappa21} and $\sup_{t \in I}\kappa_{21}(t) < \infty$ hold. Let $\varepsilon>0$ be given, and let $\bm{\phi}(t)$ satisfy \eqref{inequ}. Define $\bm{f}(t)$ by \eqref{deff}. Then $\sup_{t \in I}\|\bm{f}(t)\|_{\infty} \le \varepsilon$ holds, and $\bm{\phi}(t)$ is a solution of \eqref{perturbed}, and so that $\bm{\phi}(t)$ is written by \eqref{variation}, where $X(t)$ and $X^{-1}(t)$ are given in \eqref{funda2}. 
Now we consider a solution $\bm{x}(t)$ of \eqref{main2} defined by
\begin{equation}
 \bm{x}(t) := X(t)\bm{x}_{21}, \quad \bm{x}_{21} := X^{-1}(t_0)\bm{\phi}(t_0)+\int_{t_0}^{b}X^{-1}(s)\bm{f}(s)ds.
 \label{x_21}
\end{equation}
It will be shown later in the calculation that the (improper) integral in \eqref{x_21} converges. By \eqref{variation} and \eqref{x_21}, we obtain \eqref{phiminus11}. Let $\bm{f}(t) = (f_1(t), f_2(t))^{T}$. Then using \eqref{phiminus11}, \eqref{funda2}, and \eqref{kappa21}, we have
\begin{align*}
 \|\bm{\phi}(t)-\bm{x}(t)\|_{\infty}
  &\le  \int_t^{b}\left\|X(t)X^{-1}(s)\bm{f}(s)\right\|_{\infty}ds\\
  &= \int_t^{b}e^{-\int_t^s\Re(\lambda(\tau))d\tau}\left\|
\begin{pmatrix}
1 & -\int_t^s\mu(\tau)d\tau\\
0 & 1
\end{pmatrix}
\begin{pmatrix}
f_1(s)\\
f_2(s)
\end{pmatrix}
\right\|_{\infty}ds \\
  &= \int_t^{b}e^{-\int_t^s\Re(\lambda(\tau))d\tau}\left\|
\begin{pmatrix}
f_1(s)-f_2(s)\int_t^s\mu(\tau)d\tau\\
f_2(s)
\end{pmatrix}
\right\|_{\infty}ds \\
  &\le \int_t^{b}e^{-\int_t^s\Re(\lambda(\tau))d\tau}\max\left\{|f_1(s)|+|f_2(s)|\left|\int_t^s\mu(\tau)d\tau\right|, |f_2(s)|\right\}ds \\
  &\le \int_t^{b} \|\bm{f}(s)\|_{\infty}e^{-\int_t^s\Re(\lambda(\tau))d\tau}\left(1+\left|\int_t^s\mu(\tau)d\tau\right|\right)ds \\
  &\le \varepsilon \kappa_{21}(t)
\end{align*}
for all $t \in I$. Hence $\bm{x}_{21}$ is well-defined, and $\sup_{t \in I}\|\bm{\phi}(t)-\bm{x}(t)\|_{\infty} \le K_{21}\varepsilon$. 
Therefore, \eqref{main2} is Ulam stable on $I$, with an Ulam constant $K_{21}$. Moreover, from \eqref{variation} and \eqref{x_21}, we have
\[ \bm{x}(t) = X(t)\left(\lim_{t\to b-0} X^{-1}(t)\bm{\phi}(t)\right). \]

Next, we consider the solution to \eqref{main2} given by $\bm{y}(t) := X(t)\bm{y}_{21}$ with $\bm{y}_{21} \ne \bm{x}_{21}$. If $\sup_{t \in I}\|\bm{\phi}(t)-\bm{y}(t)\|_{\infty} \le K_{21}\varepsilon$, then
\[ \|X(t)(\bm{y}_{21}-\bm{x}_{21})\|_{\infty} \le \|\bm{\phi}(t)-\bm{y}(t)\|_{\infty} + \|\bm{\phi}(t)-\bm{x}(t)\|_{\infty} \le 2K_{21} \varepsilon \]
for all $t \in I$. However, by \eqref{unbound21}, we see that
\begin{align*}
 \lim_{t\to b-0}\|X(t)(\bm{y}_{21}-\bm{x}_{21})\|_{\infty}
  &= \lim_{t\to b-0}e^{\int_{t_0}^t\Re(\lambda(s))ds}\left\|
\begin{pmatrix}
y_1-y_2\int_{t_0}^t\mu(s)ds\\
y_2
\end{pmatrix}
\right\|_{\infty} \\
  &= \lim_{t\to b-0}e^{\int_{t_0}^t\Re(\lambda(s))ds}\max\left\{\left|y_1-y_2\int_{t_0}^t\mu(s)ds\right|, |y_2|\right\} \\
  &\ge \lim_{t\to b-0}e^{\int_{t_0}^t\Re(\lambda(s))ds} \begin{cases} |y_1| & \text{if} \quad y_2 = 0 \\ |y_2| & \text{if} \quad y_2 \ne 0 \end{cases}
  = \infty,
\end{align*}
where $\bm{y}_{21}-\bm{x}_{21} = (y_1,y_2)^T$ and $(y_1,y_2) \ne (0,0)$. This is a contradiction. Hence $\bm{x}(t)$ is the unique solution of \eqref{main2} satisfying $\sup_{t \in I}\|\bm{\phi}(t)-\bm{x}(t)\|_{\infty} \le K_{21}\varepsilon$. 

Case (ii) can be proved by combining the proofs of Case (i) and Theorem \ref{UlamS1} (ii), so we omit the proof. Thus, the proof is now complete.
\end{proof}

The lower bound of the Ulam constants are as follows.


\begin{theorem}\label{lower2}
Let $I$ be either $(a,b)$, $(a,b]$, $[a,b)$ or $[a,b]$, where $a \le a < b \le \infty$. Suppose that $\mu(t)$ is nonnegative on $I$. Then the following (i) and (ii) below hold:
\begin{itemize}
  \item[(i)] if \eqref{kappa21}, \eqref{unbound21}, and $\sup_{t \in I}\kappa_{21}(t) < \infty$ hold, then \eqref{main2} is Ulam stable on $I$, and any Ulam constant is greater than or equal to $K_{21} = \sup_{t \in I}\kappa_{21}(t)$;
  \item[(ii)] if \eqref{kappa22}, \eqref{unbound22}, and $\sup_{t \in I}\kappa_{22}(t) < \infty$ hold, then \eqref{main2} is Ulam stable on $I$, and any Ulam constant is greater than or equal to $K_{22} = \sup_{t \in I}\kappa_{22}(t)$.
\end{itemize}
\end{theorem}

\begin{proof}
Case (i). Suppose that \eqref{kappa21}, \eqref{unbound21}, and $\sup_{t \in I}\kappa_{21}(t) < \infty$ hold. Define
\begin{equation}
 \bm{f}_2(t) := \varepsilon e^{i\int_{t_0}^t\Im(\lambda(s))ds}\begin{pmatrix}
1 \\
-1
\end{pmatrix},\quad t \in I.
 \label{f2}
\end{equation}
Now we consider the solution $\bm{\phi}_1(t)$ of the system $\bm{\phi}_1' - A(t)\bm{\phi}_1 = \bm{f}_2(t)$. Since $\|\bm{f}_2(t)\|_{\infty} = \varepsilon$ holds for all $t \in I$, $\bm{\phi}_1(t)$ satisfies \eqref{inequ}. By Theorem~\ref{UlamS2}~(i), there exists the unique solution of \eqref{main2} denoted by
\[ \bm{x}_1(t) = X(t)\left(\lim_{t\to b-0} X^{-1}(t)\bm{\phi}_1(t)\right) \]
such that $\sup_{t \in I}\|\bm{\phi}_1(t)-\bm{x}_1(t)\|_{\infty} \le K_{21}\varepsilon$, where $K_{21} = \sup_{t \in I}\kappa_{21}(t)$. 
Hence, we obtain \eqref{phiminus11}. Since $\mu(t)$ is nonnegative on $I$, we see that
\begin{align*}
 \|&\bm{\phi}_1(t)-\bm{x}_1(t)\|_{\infty} = \left\|X(t)\int_t^{b}X^{-1}(s)\bm{f}_2(s)ds\right\|_{\infty} \\
  &= \varepsilon \left\|\int_t^{b} e^{-\int_t^s\lambda(\tau)d\tau+i\int_{t_0}^s\Im(\lambda(\tau))d\tau}
\begin{pmatrix}
1+\int_t^s\mu(\tau)d\tau\\
-1
\end{pmatrix}ds\right\|_{\infty}\\
  &= \varepsilon e^{-\int_t^{t_0} \Re(\lambda(\tau))d\tau} \left\|\int_t^{b} e^{-\int_{t_0}^s \Re(\lambda(\tau))d\tau}
\begin{pmatrix}
1+\int_t^s\mu(\tau)d\tau\\
-1
\end{pmatrix}ds\right\|_{\infty}\\
  &= \varepsilon e^{-\int_t^{t_0} \Re(\lambda(\tau))d\tau} \max\left\{\left|\int_t^{b} e^{-\int_{t_0}^s\Re(\lambda(\tau))d\tau}\left(1+\int_t^s\mu(\tau)d\tau\right)ds\right|, 
 \left|-\int_t^{b} e^{-\int_{t_0}^s\Re(\lambda(\tau))d\tau}ds\right|\right\}\\
  &= \varepsilon \int_t^{b} e^{-\int_t^s\Re(\lambda(\tau))d\tau}\left(1+\int_t^s\mu(\tau)d\tau\right)ds = \varepsilon \kappa_{21}(t),
\end{align*}
and so that $\sup_{t \in I}\left\|\bm{\phi}_1(t)-\bm{x}_1(t)\right\|_{\infty} = K_{21}\varepsilon$. Consequently, the Ulam constant is at least $K_{21}$. 

Case (ii) can be proved in the same way as Case (i), so we omit the proof. The proof is now complete.
\end{proof}

Using Theorems \ref{UlamS2} and \ref{lower2}, we obtain the following result, immediately. 


\begin{theorem}\label{best2}
Let $I$ be either $(a,b)$, $(a,b]$, $[a,b)$ or $[a,b]$, where $a \le a < b \le \infty$. Suppose that $\mu(t)$ is nonnegative on $I$. Then the following (i) and (ii) below hold:
\begin{itemize}
  \item[(i)] if \eqref{kappa21}, \eqref{unbound21}, and $\sup_{t \in I}\kappa_{21}(t) < \infty$ hold, then \eqref{main2} is Ulam stable on $I$, and the best Ulam constant is $K_{21} = \sup_{t \in I}\kappa_{21}(t)$;
  \item[(ii)] if \eqref{kappa22}, \eqref{unbound22}, and $\sup_{t \in I}\kappa_{22}(t) < \infty$ hold, then \eqref{main2} is Ulam stable on $I$, and the best Ulam constant is $K_{22} = \sup_{t \in I}\kappa_{22}(t)$.
\end{itemize}
\end{theorem}

When $\lambda(t) \equiv \lambda$ and $\mu(t) \equiv 1$, \eqref{main2} reduces to the constant coefficients two-dimensional linear differential system
\begin{equation}
 \bm{x}' = A\bm{x}, \quad 
A = 
\begin{pmatrix}
\lambda & 1\\
0 & \lambda
\end{pmatrix}.
 \label{const2}
\end{equation}
For this system, we have the following result.


\begin{corollary}\label{constbest2}
Let $I=\R$ and $\Re(\lambda) \ne 0$. Then \eqref{const2} is Ulam stable on $\R$, and the best Ulam constant is $K_{c2} := \frac{|\Re(\lambda)|+1}{(\Re(\lambda))^2}$. 
\end{corollary}

\begin{proof}
First, we consider the case $\Re(\lambda) > 0$. Then
\[ \kappa_{21}(t) = \int_t^{\infty} (1+s-t)e^{-\Re(\lambda)(s-t)}ds = \frac{1}{\Re(\lambda)}+\frac{1}{(\Re(\lambda))^2} = \frac{|\Re(\lambda)|+1}{(\Re(\lambda))^2}. \]
for all $t \in \R$. By Theorem~\ref{best2} (i) with $I=\R$, \eqref{const2} is Ulam stable on $\R$, and the best Ulam constant is $K_{c2}$. 

Next, we consider the case $\Re(\lambda) <0$. Then
\[ \kappa_{22}(t) = \int_{-\infty}^t (1+t-s)e^{\Re(\lambda)(t-s)}ds = \frac{1}{-\Re(\lambda)}+\frac{1}{(\Re(\lambda))^2} = \frac{|\Re(\lambda)|+1}{(\Re(\lambda))^2}. \]
for all $t \in \R$. By Theorem~\ref{best2} (ii) with $I=\R$, \eqref{const2} is Ulam stable on $\R$, and the best Ulam constant is $K_{c2}$. 
\end{proof}

\begin{remark}
When $\lambda \in \R\setminus\{0\}$, we can check that the Ulam constant
\[ \frac{|\Re(\lambda)|+1}{(\Re(\lambda))^2} = \frac{|\lambda|+1}{\lambda^2} = \|A^{-1}\|_{\infty} \]
in Corollary \ref{constbest2} is the best Ulam constant, by using the result in \cite[Theorem 4.3.]{AndOni5}.
\end{remark}


\section{Ulam stability of generalized Jordan normal form (III)}

In this section, we consider the two-dimensional linear differential system
\begin{equation}
 \bm{x}' = A(t)\bm{x}, \quad 
A(t) = 
\begin{pmatrix}
\alpha(t) & \beta(t)\\
-\beta(t) & \alpha(t)
\end{pmatrix},
 \label{main3}
\end{equation}
where $\alpha$, $\beta: I \to \R$ are continuous, and $\bm{x}\in \R^2$. A fundamental matrix of \eqref{main3} and its inverse matrix are as follows:
\begin{equation}
\begin{array}{l}
X(t) = e^{\int_{t_0}^t\alpha(s)ds}\begin{pmatrix}
\cos \int_{t_0}^t\beta(s)ds & \sin \int_{t_0}^t\beta(s)ds\\
-\sin \int_{t_0}^t\beta(s)ds & \cos \int_{t_0}^t\beta(s)ds
\end{pmatrix}, \\[5mm]
X^{-1}(t) = e^{-\int_{t_0}^t\alpha(s)ds}\begin{pmatrix}
\cos \int_{t_0}^t\beta(s)ds & -\sin \int_{t_0}^t\beta(s)ds\\
\sin \int_{t_0}^t\beta(s)ds & \cos \int_{t_0}^t\beta(s)ds
\end{pmatrix}
\end{array}
 \label{funda3}
\end{equation}
for $t \in I$, where $t_0 \in I$. 
In this section, we use the norm of the vector $\bm{v} = (v_1,v_2)^T$ with the Euclidean norm $\|\bm{v}\|_{2}:= \sqrt{|v_1|^2+|v_2|^2}$. Differences in norms do not affect Ulam stability, but they can affect Ulam constants. By referring to the proof of the theorem below, the reader will understand that the Euclidean norm is one of the appropriate norms for \eqref{main3}.


\begin{theorem}\label{UlamS3}
Let $I$ be either $(a,b)$, $(a,b]$, $[a,b)$ or $[a,b]$, where $-\infty \le a < b \le \infty$. Then the following (i) and (ii) below hold:
\begin{itemize}
  \item[(i)] if
\begin{equation}
 \kappa_{31}(t):= \int_t^{b} e^{-\int_t^s\alpha(\tau)d\tau}ds \;\; \text{exists for all}\;\; t \in I,
 \label{kappa31}
\end{equation}
and $\sup_{t \in I}\kappa_{31}(t) < \infty$, then \eqref{main3} is Ulam stable on $I$, with an Ulam constant $K_{31} := \sup_{t \in I}\kappa_{31}(t)$. 
Furthermore, if
\begin{equation}
 \lim_{t\to b-0}\int_{t_0}^t\alpha(s)ds = \infty, \quad t_0 \in (a,b),
 \label{unbound31}
\end{equation}
then for any $\varepsilon>0$ and for any continuously differentiable function $\bm{\phi}(t)$ satisfying \eqref{inequ}, there exists the unique solution of \eqref{main3} denoted by
\[ \bm{x}(t) = X(t)\left(\lim_{t\to b-0} X^{-1}(t)\bm{\phi}(t)\right) \]
such that $\sup_{t \in I}\|\bm{\phi}(t)-\bm{x}(t)\|_{2} \le K_{31}\varepsilon$, where $X(t)$ and $X^{-1}(t)$ are given in \eqref{funda3}, and $\left(\lim_{t\to b-0} X^{-1}(t)\bm{\phi}(t)\right)$ is a well-defined constant vector;
  \item[(ii)] if
\begin{equation}
 \kappa_{32}(t):= \int_{a}^t e^{\int_s^t\alpha(\tau)d\tau}ds \;\; \text{exists for all}\;\; t \in I,
 \label{kappa32}
\end{equation}
and $\sup_{t \in I}\kappa_{32}(t) < \infty$, then \eqref{main3} is Ulam stable on $I$, with an Ulam constant $K_{32} := \sup_{t \in I}\kappa_{32}(t)$. 
Furthermore, if
\begin{equation}
 \lim_{t\to a+0}\int_t^{t_0}\alpha(s)ds = -\infty, \quad t_0 \in (a,b),
 \label{unbound32}
\end{equation}
then for any $\varepsilon>0$ and for any continuously differentiable function $\bm{\phi}(t)$ satisfying \eqref{inequ}, there exists the unique solution of \eqref{main3} denoted by
\[ \bm{x}(t) = X(t)\left(\lim_{t\to a+0} X^{-1}(t)\bm{\phi}(t)\right) \]
such that $\sup_{t \in I}\|\bm{\phi}(t)-\bm{x}(t)\|_{2} \le K_{32}\varepsilon$, where $X(t)$ and $X^{-1}(t)$ are given in \eqref{funda3}, and $\left(\lim_{t\to a+0} X^{-1}(t)\bm{\phi}(t)\right)$ is a well-defined constant vector.
\end{itemize}
\end{theorem}

\begin{proof}
Since the proof of this theorem can be proved in the same way as Theorem \ref{UlamS1}, only the essential parts will be written below. Note that the maximum norm that has been used in the proofs of the previous theorems can be changed to the Euclidean norm except for computations using the definition of the maximum norm, and this proof uses everything with the Euclidean norm. 

Case (i). Suppose that \eqref{kappa31} and $\sup_{t \in I}\kappa_{31}(t) < \infty$ hold. Let $\varepsilon>0$ be given, and let $\bm{\phi}(t)$ satisfy \eqref{inequ}. Define $\bm{f}(t)$ by \eqref{deff}. Then $\sup_{t \in I}\|\bm{f}(t)\|_{2} \le \varepsilon$ holds, and $\bm{\phi}(t)$ is a solution of \eqref{perturbed}, and so that $\bm{\phi}(t)$ is written by \eqref{variation}, where $X(t)$ and $X^{-1}(t)$ are given in \eqref{funda3}. 
Now we consider a solution $\bm{x}(t)$ of \eqref{main3} defined by
\begin{equation}
 \bm{x}(t) := X(t)\bm{x}_{31}, \quad \bm{x}_{31} := X^{-1}(t_0)\bm{\phi}(t_0)+\int_{t_0}^{b}X^{-1}(s)\bm{f}(s)ds.
 \label{x_31}
\end{equation}
It will be shown later in the calculation that the (improper) integral in \eqref{x_31} converges. By \eqref{variation} and \eqref{x_31}, we obtain \eqref{phiminus11}. Let $\bm{f}(t) = (f_1(t), f_2(t))^{T}$. Then using \eqref{phiminus11}, \eqref{funda3}, \eqref{kappa31}, and the definition of the Euclidean norm, we have
\begin{align}
 \|\bm{\phi}(t)-\bm{x}(t)\|_{2}
  &\le \int_t^{b}\left\|X(t)X^{-1}(s)\bm{f}(s)\right\|_{2}ds \nonumber\\
  &= \int_t^{b}e^{-\int_t^s\alpha(\tau)d\tau}\left\|
\begin{pmatrix}
f_1(s)\cos \int_s^t\beta(\tau)d\tau-f_2(s)\sin \int_s^t\beta(\tau)d\tau\\
f_1(s)\sin \int_s^t\beta(\tau)d\tau+f_2(s)\cos \int_s^t\beta(\tau)d\tau\\
\end{pmatrix}
\right\|_{2}ds \label{norm} \\
  &= \int_t^{b}e^{-\int_t^s\alpha(\tau)d\tau}\|\bm{f}(s)\|_{2}ds \nonumber\\
  &\le \varepsilon \kappa_{31}(t) \nonumber
\end{align}
for all $t \in I$. Hence $\bm{x}_{31}$ is well-defined, and $\sup_{t \in I}\|\bm{\phi}(t)-\bm{x}(t)\|_{2} \le K_{31}\varepsilon$. 
Therefore, \eqref{main3} is Ulam stable on $I$, with an Ulam constant $K_{31}$. Moreover, from \eqref{variation} and \eqref{x_31}, we have
\[ \bm{x}(t) = X(t)\left(\lim_{t\to b-0} X^{-1}(t)\bm{\phi}(t)\right). \]

Next, we consider the solution to \eqref{main3} given by $\bm{y}(t) := X(t)\bm{y}_{31}$ with $\bm{y}_{31} \ne \bm{x}_{31}$. If $\sup_{t \in I}\|\bm{\phi}(t)-\bm{y}(t)\|_{2} \le K_{31}\varepsilon$, then
\[ \|X(t)(\bm{y}_{31}-\bm{x}_{31})\|_{2} \le \|\bm{\phi}(t)-\bm{y}(t)\|_{2} + \|\bm{\phi}(t)-\bm{x}(t)\|_{2} \le 2K_{31} \varepsilon \]
for all $t \in I$. However, by \eqref{unbound31}, we see that
\begin{align*}
 \lim_{t\to b-0}\|X(t)(\bm{y}_{31}-\bm{x}_{31})\|_{2}
  &= \lim_{t\to b-0}e^{\int_{t_0}^t\alpha(s)ds}\left\|
\begin{pmatrix}
y_1\cos \int_{t_0}^t\beta(s)ds+y_2\sin \int_{t_0}^t\beta(s)ds\\
-y_1\sin \int_{t_0}^t\beta(s)ds+y_2\cos \int_{t_0}^t\beta(s)ds\\
\end{pmatrix}
\right\|_{2} \\
  &= \lim_{t\to b-0}e^{\int_{t_0}^t\alpha(s)ds}\|(\bm{y}_{31}-\bm{x}_{31})\|_{2}
  = \infty, 
\end{align*}
where $\bm{y}_{31}-\bm{x}_{31} = (y_1,y_2)^T$ and $(y_1,y_2) \ne (0,0)$. This is a contradiction. Hence $\bm{x}(t)$ is the unique solution of \eqref{main3} satisfying $\sup_{t \in I}\|\bm{\phi}(t)-\bm{x}(t)\|_{2} \le K_{31}\varepsilon$. 

Case (ii) can be proved by combining the proofs of Case (i) and Theorem \ref{UlamS1} (ii), so we omit the proof. Thus, the proof is now complete.
\end{proof}

The lower bound of the Ulam constants are as follows.


\begin{theorem}\label{lower3}
Let $I$ be either $(a,b)$, $(a,b]$, $[a,b)$ or $[a,b]$, where $a \le a < b \le \infty$. Then the following (i) and (ii) below hold:
\begin{itemize}
  \item[(i)] if \eqref{kappa31}, \eqref{unbound31}, and $\sup_{t \in I}\kappa_{31}(t) < \infty$ hold, then \eqref{main3} is Ulam stable on $I$, and any Ulam constant is greater than or equal to $K_{31} = \sup_{t \in I}\kappa_{31}(t)$;
  \item[(ii)] if \eqref{kappa32}, \eqref{unbound32}, and $\sup_{t \in I}\kappa_{32}(t) < \infty$ hold, then \eqref{main3} is Ulam stable on $I$, and any Ulam constant is greater than or equal to $K_{32} = \sup_{t \in I}\kappa_{32}(t)$.
\end{itemize}
\end{theorem}

\begin{proof}
Case (i). Suppose that \eqref{kappa31}, \eqref{unbound31}, and $\sup_{t \in I}\kappa_{31}(t) < \infty$ hold. Let $\bm{x}_* \in \R^2\setminus\{\bm{0}\}$. Define
\begin{equation}
 \bm{f}_3(t) := \frac{\varepsilon e^{-\int_{t_0}^t\alpha(s)ds}}{\|\bm{x}_*\|_{2}} X(t)\bm{x}_*,\quad t \in I,
 \label{f3}
\end{equation}
where $X(t)$ is given in \eqref{funda3}. Now we consider the solution $\bm{\phi}_1(t)$ of the system $\bm{\phi}_1' - A(t)\bm{\phi}_1 = \bm{f}_3(t)$. Since $\|\bm{f}_3(t)\|_{2} = \varepsilon$ holds for all $t \in I$, $\bm{\phi}_1(t)$ satisfies \eqref{inequ}. By Theorem~\ref{UlamS3}~(i), there exists the unique solution of \eqref{main3} denoted by
\[ \bm{x}_1(t) = X(t)\left(\lim_{t\to b-0} X^{-1}(t)\bm{\phi}_1(t)\right) \]
such that $\sup_{t \in I}\|\bm{\phi}_1(t)-\bm{x}_1(t)\|_{2} \le K_{31}\varepsilon$, where $K_{31} = \sup_{t \in I}\kappa_{31}(t)$. 
Hence, we obtain \eqref{phiminus11}. Since $\left\|X(t)\bm{x}_*\right\|_{2}=e^{\int_{t_0}^t\alpha(s)ds}\|\bm{x}_*\|_{2}$ holds, we see that
\begin{align*}
 \|\bm{\phi}_1(t)-\bm{x}_1(t)\|_{2} &= \left\|X(t)\int_t^{b}X^{-1}(s)\bm{f}_3(s)ds\right\|_{2} = \left\|X(t)\int_t^{b}\frac{\varepsilon e^{-\int_{t_0}^s\alpha(\tau)d\tau}}{\|\bm{x}_*\|_{2}} \bm{x}_*ds\right\|_{2} \\
  &= \left\|X(t)\bm{x}_*\right\|_{2}\int_t^{b}\frac{\varepsilon e^{-\int_{t_0}^s\alpha(\tau)d\tau}}{\|\bm{x}_*\|_{2}} ds
  = \varepsilon \kappa_{31}(t),
\end{align*}
and so that $\sup_{t \in I}\left\|\bm{\phi}_1(t)-\bm{x}_1(t)\right\|_{2} = K_{31}\varepsilon$. Consequently, the Ulam constant is at least $K_{31}$. 

Case (ii) can be proved in the same way as Case (i), so we omit the proof. The proof is now complete.
\end{proof}

Using Theorems \ref{UlamS3} and \ref{lower3}, we obtain the following result, immediately. 


\begin{theorem}\label{best3}
Let $I$ be either $(a,b)$, $(a,b]$, $[a,b)$ or $[a,b]$, where $a \le a < b \le \infty$. Then the following (i) and (ii) below hold:
\begin{itemize}
  \item[(i)] if \eqref{kappa31}, \eqref{unbound31}, and $\sup_{t \in I}\kappa_{31}(t) < \infty$ hold, then \eqref{main3} is Ulam stable on $I$, and the best Ulam constant is $K_{31} = \sup_{t \in I}\kappa_{31}(t)$;
  \item[(ii)] if \eqref{kappa32}, \eqref{unbound32}, and $\sup_{t \in I}\kappa_{32}(t) < \infty$ hold, then \eqref{main3} is Ulam stable on $I$, and the best Ulam constant is $K_{32} = \sup_{t \in I}\kappa_{32}(t)$.
\end{itemize}
\end{theorem}

When $\alpha(t) \equiv \alpha$ and $\beta(t) \equiv \beta$, \eqref{main3} reduces to the constant coefficients two-dimensional linear differential system
\begin{equation}
 \bm{x}' = A\bm{x}, \quad 
A = 
\begin{pmatrix}
\alpha & \beta\\
-\beta & \alpha
\end{pmatrix}.
 \label{const3}
\end{equation}
For this system, we have the following result.


\begin{corollary}\label{constbest3}
Let $I=\R$ and $\alpha \ne 0$. Then \eqref{const3} is Ulam stable on $\R$, and the best Ulam constant is $K_{c3} := \frac{1}{|\alpha|}$. 
\end{corollary}

\begin{proof}
First, we consider the case $\alpha > 0$. Then
\[ \kappa_{31}(t) = \int_t^{\infty} e^{-\alpha(s-t)}ds = \frac{1}{\alpha} = \frac{1}{|\alpha|}, \]
for all $t \in \R$. By Theorem~\ref{best3} (i) with $I=\R$, \eqref{const3} is Ulam stable on $\R$, and the best Ulam constant is $K_{c3}$. 

Next, we consider the case $\alpha <0$. Then
\[ \kappa_{32}(t) = \int_{-\infty}^t e^{\alpha(t-s)}ds = \frac{1}{-\alpha} = \frac{1}{|\alpha|}, \]
for all $t \in \R$. By Theorem~\ref{best3} (ii) with $I=\R$, \eqref{const3} is Ulam stable on $\R$, and the best Ulam constant is $K_{c3}$. 
\end{proof}


\begin{remark}
If we use the maximum norm instead of the Euclidean norm in \eqref{norm} in the proof of Theorem \ref{UlamS3}, then we obtain
\begin{align*}
 \|\bm{\phi}(t)-\bm{x}(t)\|_{\infty}
  &\le \int_t^{b}e^{-\int_t^s\alpha(\tau)d\tau}\left\|
\begin{pmatrix}
f_1(s)\cos \int_s^t\beta(\tau)d\tau-f_2(s)\sin \int_s^t\beta(\tau)d\tau\\
f_1(s)\sin \int_s^t\beta(\tau)d\tau+f_2(s)\cos \int_s^t\beta(\tau)d\tau\\
\end{pmatrix}
\right\|_{\infty}ds \\
  &= \int_t^{b}e^{-\int_t^s\alpha(\tau)d\tau}\sqrt{(f_1(s))^2+(f_2(s))^2}\\
  &\hspace{5mm} \times \max\left\{\left|\sin \left(\frac{\pi}{2}+\int_s^t\beta(\tau)d\tau+\gamma(\tau)\right)\right|, \left|\sin \left(\int_s^t\beta(\tau)d\tau+\gamma(\tau)\right)\right|\right\}ds \\
  &\le \varepsilon \sqrt{2}\kappa_{31}(t),
\end{align*}
where
\[ \frac{f_2(s)}{\sqrt{(f_1(s))^2+(f_2(s))^2}} = \sin \gamma(t), \quad \frac{f_1(s)}{\sqrt{(f_1(s))^2+(f_2(s))^2}} = \cos \gamma(t) \]
for $t\in I$. This says that an Ulam constant is $\sqrt{2}K_{31}$, when we use the maximum norm. On the other hand, if we use the maximum norm instead of the Euclidean norm in \eqref{f3} in the proof of Theorem \ref{lower3}, and we choose $\bm{x}_*=(1,-1)^T$, then we obtain
\begin{align*}
 \|&\bm{\phi}_1(t)-\bm{x}_1(t)\|_{\infty} = \left\|X(t)\bm{x}_*\right\|_{\infty}\int_t^{b}\frac{\varepsilon e^{-\int_{t_0}^s\alpha(\tau)d\tau}}{\|\bm{x}_*\|_{\infty}} ds\\
  &= \left\|\begin{pmatrix}
\cos \int_{t_0}^t\beta(s)ds - \sin \int_{t_0}^t\beta(s)ds\\
-\sin \int_{t_0}^t\beta(s)ds - \cos \int_{t_0}^t\beta(s)ds
\end{pmatrix}\right\|_{\infty}\int_t^{b}\varepsilon e^{-\int_t^s\alpha(\tau)d\tau} ds \\
  &= \sqrt{2}\max\left\{\left|\sin \left(\int_{t_0}^t\beta(s)ds+\frac{3\pi}{4}\right)\right|, \left|\sin \left(\int_{t_0}^t\beta(s)ds+\frac{\pi}{4}\right)\right|\right\}\int_t^{b}\varepsilon e^{-\int_t^s\alpha(\tau)d\tau} ds
\end{align*}
for $t\in I$, and thus, $\|\bm{\phi}_1(t)-\bm{x}_1(t)\|_{\infty} \ge \varepsilon \kappa_{31}(t)$ for $t\in I$. Moreover, we assume that $I=\R$, $\alpha(t) \equiv \alpha > 0$ and $\beta(t) \equiv \beta \ne 0$. Then we see that
\[ \|\bm{\phi}_1(t)-\bm{x}_1(t)\|_{\infty} = \frac{\sqrt{2}\varepsilon}{\alpha}\max\left\{\left|\sin \left(\beta(t-t_0)+\frac{3\pi}{4}\right)\right|, \left|\sin \left(\beta(t-t_0)+\frac{\pi}{4}\right)\right|\right\} \]
for $t\in I$. Define $t_n=\frac{1}{\beta}\left(\frac{\pi}{4}+n\pi\right)+t_0$ for $n \in \Z$. If we choose $t=t_n$ for $n \in \Z$, we have
\[ \|\bm{\phi}_1(t_n)-\bm{x}_1(t_n)\|_{\infty} = \frac{\sqrt{2}\varepsilon}{\alpha}, \quad n \in \Z. \]
Therefore, any Ulam constant is greater than or equal to $\frac{\sqrt{2}}{\alpha}$, and so that the best Ulam constant for \eqref{const3} when using the maximum norm is $\frac{\sqrt{2}}{\alpha}$. 
\end{remark}

By the above remark, the following result holds.


\begin{theorem}\label{constbest4}
Let $I=\R$ and $\alpha \ne 0$ and $\beta \ne 0$. Then \eqref{const3} is Ulam stable on $\R$, and the best Ulam constant when using the maximum norm is $K_{c4} := \frac{\sqrt{2}}{|\alpha|}$. 
\end{theorem}

\begin{proof}
The proof for $\alpha < 0$ is similar to that for $\alpha > 0$. Therefore, the proof is omitted. 
\end{proof}


\begin{remark}
Note that Theorem \ref{constbest4} is also a new theorem. In \cite{AndOni5}, the lower bound of the Ulam constant for \eqref{const3} was not proved. In addition, from the conclusions of Corollary~\ref{constbest3} and Theorem~\ref{constbest4}, we find that the best Ulam constants vary depending on the choice of norm. 
\end{remark}


\section{Examples and approximations of node, saddle, and focus}

\begin{example}\label{ex:blowup}
In \eqref{main1}, let $I=(0, 1)$, $A(t) = 
\begin{pmatrix}
\frac{1}{1-t} & 0 \\
0 & -\frac{1}{t}
\end{pmatrix}$, and $t_0 \in I$. Then
$$ X(t)= \begin{pmatrix} \frac{1-t_0}{1-t} & 0 \\ 0 & \frac{t_0}{t} \end{pmatrix} \quad\text{and}\quad X^{-1}(t) = \begin{pmatrix} \frac{1-t}{1-t_0} & 0 \\ 0 & \frac{t}{t_0} \end{pmatrix}, $$
and all conditions of Theorem \ref{UlamS1} (iii) are satisfied, namely:
\begin{itemize}
  \item[(i)] Since
\[ \kappa_{13}(t) := \max\left\{\int_t^1 e^{-\int_t^s\frac{1}{1-\tau}d\tau}ds, \int_0^t e^{\int_s^t\left(-\frac{1}{\tau}\right)d\tau}ds\right\} 
= \frac{1}{2}\left(\left|t-\frac{1}{2}\right|+\frac{1}{2}\right) \]
exists for all $t \in I$, and $\sup_{t \in I}\kappa_{13}(t) = \frac{1}{2} < \infty$, then \eqref{main1} is Ulam stable on $I$, with an Ulam constant $K_{13}:=\frac{1}{2}$. Furthermore, since
\[ \lim_{t\to 1-0}\int_{t_0}^t\frac{1}{1-s}ds = \lim_{t\to 1-0}\log\frac{1-t_0}{1-t}= \infty, \quad \lim_{t\to +0}\int_{t_0}^t\frac{1}{s}ds = \lim_{t\to +0}\log\frac{t_0}{t}= \infty, \]
then for any $\varepsilon>0$ and for any continuously differentiable function $\bm{\phi}(t) = (\phi_1(t),\phi_2(t))^T$ satisfying \eqref{inequ}, there exists the unique solution of \eqref{main1} denoted by
\[ \bm{x}(t) = X(t)\begin{pmatrix}
\lim_{t\to 1-0} \frac{(1-t)\phi_1(t)}{1-t_0}\\
\lim_{t\to +0} \frac{t\phi_2(t)}{t_0}
\end{pmatrix} \]
such that $\sup_{t \in I}\|\bm{\phi}(t)-\bm{x}(t)\|_{\infty} \le K_{13}\varepsilon$, where $\begin{pmatrix}
\lim_{t\to 1-0} \frac{(1-t)\phi_1(t)}{1-t_0}\\
\lim_{t\to +0} \frac{t\phi_2(t)}{t_0}
\end{pmatrix}$ is a well-defined constant vector. 
\end{itemize}
In addition, by Theorem \ref{best1} (iii), $K_{13}=\frac{1}{2}$ is the best Ulam constant for \eqref{main1} on $I$.
\end{example}


\begin{remark}
Example \ref{ex:blowup} illustrates that our results (Theorems \ref{UlamS1}--\ref{best1}, \ref{UlamS2}--\ref{best2}, and \ref{UlamS3}--\ref{best3}) also apply to systems that have blow-up solutions. To the best of the author's knowledge, no studies of Ulam stability for equations with blow-up solutions are known so far.
\end{remark}


\begin{remark}
In many results \cite{BacDra1,BacDra2,BacDraPitSin,BusLupO'R,BusO'RSaiTab}, Ulam stability is established under the assumption that the linear part (in nonlinear system) admits an exponential dichotomy on $I=[0,\infty)$. For the discrete problem it corresponds to a discrete dichotomy, for periodic coefficients case it corresponds to that the monodromy matrix is hyperbolic. However, Example \ref{ex:blowup} suggests that the exponential dichotomy is not required for our results. 
\end{remark}

\begin{example}
In \eqref{main2}, let $I=\R$ and 
$$ A(t) = A_{1}(t)= 
\begin{pmatrix}
   1+i\beta(t) & \frac{2}{\sqrt{\pi}}e^{-t^2} \\ 0 & 1+i\beta(t)
\end{pmatrix} \quad\text{or}\quad A(t) = A_{-1}(t)= 
\begin{pmatrix}
  -1+i\beta(t) & \frac{2}{\sqrt{\pi}}e^{-t^2} \\ 0 & -1+i\beta(t)
\end{pmatrix} $$
for $t\in I$, where $\beta:I\rightarrow\R$ is continuous. Note that either
\begin{align}
 X(t)=X_{1}(t) &= e^{t}e^{i\int_0^t\beta(s)ds} \begin{pmatrix} 1 & \operatorname{erf}(t) \\ 0 & 1 \end{pmatrix}, \quad X_{1}^{-1}(t)=e^{-t}e^{-i\int_0^t\beta(s)ds} \begin{pmatrix} 1 & -\operatorname{erf}(t) \\ 0 & 1 \end{pmatrix}, \quad\text{or} \label{Xplus1def} \\
 X(t)=X_{-1}(t) &= e^{-t}e^{i\int_0^t\beta(s)ds} \begin{pmatrix} 1 & \operatorname{erf}(t) \\ 0 & 1 \end{pmatrix}, \quad X_{-1}^{-1}(t)=e^{t}e^{-i\int_0^t\beta(s)ds} \begin{pmatrix} 1 & -\operatorname{erf}(t) \\ 0 & 1 \end{pmatrix}, \label{Xminus1def}
\end{align}
where $\operatorname{erf}$ is the error function. Then Theorem \ref{UlamS2} (i) holds for $A_{1}$ and (ii) holds for $A_{-1}$. In particular, we have the following details.
\begin{itemize}
  \item[(i)] Since
\begin{equation}
 \kappa_{21}(t):= \int_t^{\infty} \left(1+\left|\int_t^s\frac{2}{\sqrt{\pi}}e^{-\tau^2}d\tau\right|\right)e^{-\int_t^s(1)d\tau}ds =  1+e^{\frac{1}{4}+t}\operatorname{erfc}\left(\frac{1}{2}+t\right)
 \label{kappa21-examp}
\end{equation}
exists for all $t \in I$, where $\operatorname{erfc}$ is the complementary error function, 
and $\sup_{t \in I}\kappa_{21}(t) = 1.78395 < \infty$ at $t = -0.603489$, then \eqref{main2} is Ulam stable on $I$, with an Ulam constant $K_{21}:=1.78395$. 
Furthermore, since
\begin{equation}
 \lim_{t\to \infty}\int_{t_0}^t (1)ds = \infty, \quad t_0 \in (-\infty,\infty),
 \label{unbound21-examp}
\end{equation}
then for any $\varepsilon>0$ and for any continuously differentiable function $\bm{\phi}(t)$ satisfying \eqref{inequ}, there exists the unique solution of \eqref{main2} denoted by
\[ \bm{x}(t) = X_1(t)\left(\lim_{t\to\infty} X_1^{-1}(t)\bm{\phi}(t)\right) \]
such that $\sup_{t \in I}\|\bm{\phi}(t)-\bm{x}(t)\|_{\infty} \le K_{21}\varepsilon=1.78395\varepsilon$, where $X_1(t)$ and $X_1^{-1}(t)$ are given in \eqref{Xplus1def}, and $\left(\lim_{t\to\infty} X_1^{-1}(t)\bm{\phi}(t)\right)$ is a well-defined constant vector.
  \item[(ii)] Since
\begin{equation}
 \kappa_{22}(t):= \int_{-\infty}^t \left(1+\left|\int_s^t\frac{2}{\sqrt{\pi}}e^{-\tau^2}d\tau\right|\right)e^{\int_s^t(-1)d\tau}ds = 1+e^{\frac{1}{4}-t}\operatorname{erfc}\left(\frac{1}{2}-t\right)
 \label{kappa22-examp}
\end{equation}
exists for all $t \in I$, and $\sup_{t \in I}\kappa_{22}(t) = 1.78395 < \infty$ at $t=0.603489$, then \eqref{main2} is Ulam stable on $I$, with an Ulam constant $K_{22}=1.78395$. 
Furthermore, since
\begin{equation}
 \lim_{t\to -\infty}\int_t^{t_0}(-1)ds = -\infty, \quad t_0\in(-\infty,\infty)
 \label{unbound22-examp}
\end{equation}
then for any $\varepsilon>0$ and for any continuously differentiable function $\bm{\phi}(t)$ satisfying \eqref{inequ}, there exists the unique solution of \eqref{main2} denoted by
\[ \bm{x}(t) = X_{-1}(t)\left(\lim_{t\to -\infty} X_{-1}^{-1}(t)\bm{\phi}(t)\right) \]
such that $\sup_{t \in I}\|\bm{\phi}(t)-\bm{x}(t)\|_{\infty} \le K_{22}\varepsilon =1.78395\varepsilon$, where $X_{-1}(t)$ and $X_{-1}^{-1}(t)$ are given in \eqref{Xminus1def}, and $\left(\lim_{t\to -\infty} X_{-1}^{-1}(t)\bm{\phi}(t)\right)$ is a well-defined constant vector.
\end{itemize}
\end{example}

\begin{example}
In \eqref{main2}, let $A(t)=\begin{pmatrix} \frac{2t}{1+t^2} & 1 \\ 0 & \frac{2t}{1+t^2} \end{pmatrix}$ and $I=\R$. Then $X(t)=(1+t^2)\begin{pmatrix} 1 & t \\ 0 & 1 \end{pmatrix}$ is a fundamental matrix for \eqref{main2}. Given $\varepsilon>0$, take $\bm{f}(t)=\begin{pmatrix} \varepsilon \\ 0 \end{pmatrix}$ and $\bm{\phi}(t)=X(t)\begin{pmatrix} \varepsilon\tan^{-1}t \\ 0 \end{pmatrix}$. It follows that $\bm{\phi}$ satisfies \eqref{deff}, and for any solution $\bm{x}(t)=X(t)\bm{x}_0$ of \eqref{main2}, we have that
$$ \|\bm{\phi}(t)-\bm{x}(t)\|_{\infty} = \left\|X(t) \left(\begin{pmatrix} \varepsilon\tan^{-1}t \\ 0 \end{pmatrix}-\bm{x}_0\right)\right\|_{\infty} $$
is unbounded on $I=(-\infty,\infty)$ for any choice of $\bm{x}_0$. Thus \eqref{main2} is not Ulam stable on $I$ for this $A$. 

Note that, since $\lambda(t) = \frac{2t}{1+t^2}$ and $\mu(t)\equiv 1$, we have
\begin{align*}
 \int_t^{T} &\left(1+\left|\int_t^s\mu(\tau)d\tau\right|\right)e^{-\int_t^s\Re(\lambda(\tau))d\tau}ds = (1+t^2)\int_t^{T} \frac{1+s-t}{1+s^2}ds \\
  &= \frac{1+t^2}{2}\left(\log\sqrt{\frac{1+T^2}{1+t^2}}+(1-t)\left(\tan^{-1}T-\tan^{-1}t\right)\right)
\end{align*}
and
\begin{align*}
 \int_{-T}^t &\left(1+\left|\int_s^t\mu(\tau)d\tau\right|\right)e^{\int_s^t\Re(\lambda(\tau))d\tau}ds = (1+t^2)\int_t^{T} \frac{1+s-t}{1+s^2}ds \\
  &= \frac{1+t^2}{2}\left(\log\sqrt{\frac{1+T^2}{1+t^2}}+(1+t)\left(\tan^{-1}t+\tan^{-1}T\right)\right)
\end{align*}
for $T>0$. When $T \to \infty$, these left-hand integrals do not converge for any $t\in I$; that is, conditions \eqref{kappa21} and \eqref{kappa22} are not satisfied. Therefore, it is concluded that \eqref{kappa21} and \eqref{kappa22} are sharp conditions.
\end{example}

\begin{example}
In \eqref{main3}, let $I=\R$, $\beta:I\rightarrow\R$ be continuous, and
$$ \alpha(t) = \alpha_1(t) = 1-\frac{2t}{1+t^2} \quad\text{or}\quad \alpha(t) = \alpha_{-1}(t) = -1-\frac{2t}{1+t^2} $$ 
for $t\in I$. Then by Theorem \ref{UlamS3}, the following (i) and (ii) below hold.
\begin{itemize}
  \item[(i)] For $\alpha=\alpha_1$, since
\begin{equation}
 \kappa_{31}(t):= \int_t^{\infty} e^{-\int_t^s\left(1-\frac{2\tau}{1+\tau^2}\right)d\tau}ds = \frac{t^2+2t+3}{t^2+1}
 \label{kappa31-examp}
\end{equation}
exists for all $t \in I$, 
and $\sup_{t \in I}\kappa_{31}(t) = 2+\sqrt{2} < \infty$ at $t=\sqrt{2}-1$, then \eqref{main3} is Ulam stable on $I$, with an Ulam constant $K_{31}:= 2+\sqrt{2}$. 
Furthermore, as
\begin{equation}
 \lim_{t\to \infty}\int_{t_0}^t\alpha(s)ds = \lim_{t\to \infty}\int_{t_0}^t \left(1-\frac{2s}{1+s^2}\right) ds = \infty, \quad t_0 \in (-\infty,\infty),
 \label{unbound31-examp}
\end{equation}
then for any $\varepsilon>0$ and for any continuously differentiable function $\bm{\phi}(t)$ satisfying \eqref{inequ}, there exists the unique solution of \eqref{main3} denoted by
\[ \bm{x}(t) = X(t)\left(\lim_{t\to\infty} X^{-1}(t)\bm{\phi}(t)\right) \]
such that $\sup_{t \in I}\|\bm{\phi}(t)-\bm{x}(t)\|_{2} \le K_{31}\varepsilon=\left(2+\sqrt{2}\right)\varepsilon$, where $X(t)$ and $X^{-1}(t)$ are given by
$$\begin{array}{l}
X(t) = \frac{e^{t-t_0}\left(1+t_0^2\right)}{1+t^2}\begin{pmatrix}
\cos \int_{t_0}^t\beta(s)ds & \sin \int_{t_0}^t\beta(s)ds\\
-\sin \int_{t_0}^t\beta(s)ds & \cos \int_{t_0}^t\beta(s)ds
\end{pmatrix}, \\[5mm]
X^{-1}(t) = \frac{e^{t_0-t}\left(1+t^2\right)}{1+t_0^2}\begin{pmatrix}
\cos \int_{t_0}^t\beta(s)ds & -\sin \int_{t_0}^t\beta(s)ds\\
\sin \int_{t_0}^t\beta(s)ds & \cos \int_{t_0}^t\beta(s)ds
\end{pmatrix}
\end{array}$$
respectively, and $\left(\lim_{t\to\infty} X^{-1}(t)\bm{\phi}(t)\right)$ is a well-defined constant vector.
  \item[(ii)] For $\alpha=\alpha_{-1}$, since
\begin{equation}
 \kappa_{32}(t):= \int_{-\infty}^t e^{\int_s^t\alpha(\tau)d\tau}ds = \frac{t^2-2t+3}{t^2+1} 
 \label{kappa32-examp}
\end{equation}
exists for all $t \in I$, and $\sup_{t \in I}\kappa_{32}(t) = 2+\sqrt{2} < \infty$ at $t=1-\sqrt{2}$, then \eqref{main3} is Ulam stable on $I$, with an Ulam constant $K_{32}:=2+\sqrt{2}$. 
Furthermore, as
\begin{equation}
 \lim_{t\to -\infty}\int_t^{t_0}\alpha(s)ds = -\infty, \quad t_0 \in (-\infty,\infty),
 \label{unbound32-examp}
\end{equation}
then for any $\varepsilon>0$ and for any continuously differentiable function $\bm{\phi}(t)$ satisfying \eqref{inequ}, there exists the unique solution of \eqref{main3} denoted by
\[ \bm{x}(t) = X(t)\left(\lim_{t\to -\infty} X^{-1}(t)\bm{\phi}(t)\right) \]
such that $\sup_{t \in I}\|\bm{\phi}(t)-\bm{x}(t)\|_{2} \le K_{32}\varepsilon = \left(2+\sqrt{2}\right)\varepsilon$, where $X(t)$ and $X^{-1}(t)$ are given by
$$\begin{array}{l}
X(t) = \frac{e^{t_0-t}\left(1+t_0^2\right)}{1+t^2}\begin{pmatrix}
\cos \int_{t_0}^t\beta(s)ds & \sin \int_{t_0}^t\beta(s)ds\\
-\sin \int_{t_0}^t\beta(s)ds & \cos \int_{t_0}^t\beta(s)ds
\end{pmatrix}, \\[5mm]
X^{-1}(t) = \frac{e^{t-t_0}\left(1+t^2\right)}{1+t_0^2}\begin{pmatrix}
\cos \int_{t_0}^t\beta(s)ds & -\sin \int_{t_0}^t\beta(s)ds\\
\sin \int_{t_0}^t\beta(s)ds & \cos \int_{t_0}^t\beta(s)ds
\end{pmatrix}
\end{array}$$
respectively, and $\left(\lim_{t\to -\infty} X^{-1}(t)\bm{\phi}(t)\right)$ is a well-defined constant vector.
\end{itemize}
\end{example}

\begin{example}\label{ex:blowup2}
Let $I=\left(0, \frac{\pi}{2}\right)$. Consider \eqref{perturbed} with
\begin{equation}
 A(t) =
\begin{pmatrix}
2 \cot t & i \\
-3i & -2 \cot t
\end{pmatrix},
 \label{nonJordan}
\end{equation}
where $\bm{f}: I \to \C^n$ is an arbitrary continuous vector function on $I$. By Theorem \ref{iff1}, we see that \eqref{perturbed} with \eqref{nonJordan} is Ulam stable on $I$, with an Ulam constant $K$ if and only if \eqref{unperturbed} with \eqref{nonJordan} is Ulam stable on $I$, and an Ulam constant is the same $K$. Therefore, we only need to consider the Ulam  stability of \eqref{unperturbed} with \eqref{nonJordan}. However, note here that $A(t)$ does not belong to the generalized Jordan normal forms (I)--(III). Let
\[ R(t) = 
\begin{pmatrix}
\cos t & i\sin t \\
i\sin t & \cos t
\end{pmatrix}. \]
Using $\bm{y}=R(t)\bm{x}$, \eqref{unperturbed} is transformed into \eqref{Jordan}. From
\begin{equation}
 J(t) = \left(R'(t)+R(t)A(t)\right)R^{-1}(t) = \begin{pmatrix}
\csc t \sec t & 0 \\
0 & -\csc t \sec t
\end{pmatrix},
 \label{Jordan-Ex}
\end{equation}
we see that \eqref{Jordan} with \eqref{Jordan-Ex} is just \eqref{main1} with $A(t)=J(t)$. Let $t_0 \in I$. Then
$$ X(t)= \begin{pmatrix} \tan t \cot t_0 & 0 \\ 0 & \cot t \tan t_0 \end{pmatrix} \quad\text{and}\quad X^{-1}(t) = \begin{pmatrix} \cot t\tan t_0 & 0 \\ 0 & \tan t\cot t_0 \end{pmatrix}, $$
and all conditions of Theorem \ref{UlamS1} (iii) are satisfied, namely:
\begin{itemize}
  \item[(i)] Since
\begin{align*}
 \kappa_{13}(t)
 &:= \max\left\{\int_t^{\frac{\pi}{2}} e^{-\int_t^s\csc \tau \sec\tau d\tau}ds, \int_0^t e^{\int_s^t\left(-\csc\tau\sec\tau\right)d\tau}ds\right\} \\
 & = \max\{-\tan t \ln(\sin t),\; -\cot t\ln(\cos t)\} \\
 & =\begin{cases} -\tan t\ln(\sin t) &: t\in\left(0,\frac{\pi}{4}\right] \\ -\cot t\ln(\cos t) &: t\in\left[\frac{\pi}{4},\frac{\pi}{2}\right) \end{cases}
\end{align*}
exists for all $t \in I$, and $\sup_{t \in I}\kappa_{13}(t) = 0.4023711 < \infty$, then \eqref{main1} with $A(t)=J(t)$ is Ulam stable on $I$, with an Ulam constant $K_{13}:=0.4023711$. Furthermore, since
\begin{align*}
\lim_{t\to \frac{\pi}{2}-0}\int_{t_0}^t\csc s\sec sds &= \lim_{t\to \frac{\pi}{2}-0} \ln(\tan t\cot t_0)= \infty, \\
 \lim_{t\to +0}\int_{t_0}^t\left(-\csc s\sec s\right)ds &= \lim_{t\to +0}\ln(\cot t\tan t_0)= \infty,
\end{align*}
then for any $\varepsilon>0$ and for any continuously differentiable function $\bm{\phi}(t) = (\phi_1(t),\phi_2(t))^T$ satisfying \eqref{inequ}, there exists the unique solution of \eqref{main1} with $A(t)=J(t)$ denoted by
\[ \bm{x}(t) = X(t)\begin{pmatrix}
\lim_{t\to \frac{\pi}{2}-0} \phi_1(t)\cot t\tan t_0\\
\lim_{t\to +0} \phi_2(t)\tan t\cot t_0
\end{pmatrix} \]
such that $\sup_{t \in I}\|\bm{\phi}(t)-\bm{x}(t)\|_{\infty} \le K_{13}\varepsilon=0.4023711 \varepsilon$, where $\begin{pmatrix}
\lim_{t\to \frac{\pi}{2}-0} \phi_1(t)\cot t\tan t_0\\
\lim_{t\to +0}  \phi_2(t)\tan t\cot t_0
\end{pmatrix}$ is a well-defined constant vector.
\end{itemize}
By Theorem \ref{best1} (iii), $K_{13}=0.4023711 $ is the best Ulam constant for \eqref{main1} with $A(t)=J(t)$ on $I$. Moreover, since
\[ \sup_{t \in I}\|R(t)\|_{\infty} = \sup_{t \in I}\|R^{-1}(t)\|_{\infty} = \sqrt{2}\sup_{t \in I}\sin\left(t+\frac{\pi}{4}\right) = \sqrt{2}, \]
$\|R(t)\|_{\infty}$ and $\|R^{-1}(t)\|_{\infty}$ are bounded on $I$. From Theorem \ref{Jimpliesu}, \eqref{unperturbed} with \eqref{nonJordan} is Ulam stable on $I$, and an Ulam constant is $2K_{13} = 0.8047422$. In addition, Theorem \ref{iff1}, \eqref{perturbed} with \eqref{nonJordan} is also Ulam stable on $I$, and an Ulam constant is $2K_{13} = 0.8047422$. 
\end{example}

Next, we propose that Ulam stability gives approximations of node, saddle, and focus. 


\begin{example}\label{ex:saddle}
In \eqref{const1}, let $I=(-\infty, \infty)$, and let
\begin{equation}
 A = 
\begin{pmatrix}
1 & 0 \\
0 & -1
\end{pmatrix}.
 \label{saddle}
\end{equation}
Then, by Corollary \ref{constbest1}, \eqref{const1} with \eqref{saddle} is Ulam stable on $I$, and the best Ulam constant is $K_{c1} = 1$. Moreover, we consider \eqref{deff} with \eqref{saddle} and
\begin{equation}
\bm{f}(t) = 
\begin{pmatrix}
0.2\left(1-2\max \{\cos t,0\}\right) \\
0
\end{pmatrix}.
 \label{perturbation}
\end{equation}
Let $\varepsilon=0.2$. Clearly the solution $\bm{\phi}(t)$ of \eqref{deff} with \eqref{saddle} and \eqref{perturbation} satisfies \eqref{inequ}. Since \eqref{kappa13}, \eqref{unbound13} and $\sup_{t \in I}\kappa_{13}(t) = K_{13} = K_{c1} = 1 < \infty$ hold, by Theorem \ref{UlamS1} (iii), we see that there exists the unique solution $\bm{x}(t)$ of \eqref{const1} with \eqref{saddle} such that $\sup_{t \in I}\|\bm{\phi}(t)-\bm{x}(t)\|_{\infty} \le K_{13}\varepsilon = 0.2$. This means that there is only one orbit of \eqref{const1} with \eqref{saddle} in the $\varepsilon$-neighborhood of any orbit of \eqref{deff} with \eqref{saddle} and \eqref{perturbation}. Since the geometric classification of the origin of \eqref{const1} with \eqref{saddle} is a saddle point (see \cite{Hale,KelPet,MicHouLiu}), a phase portrait of \eqref{deff} with \eqref{saddle} and \eqref{perturbation} can be proposed as an approximation of the saddle point (see Figure \ref{fig:1}). 
\begin{figure}[ht]
\begin{center}
\begin{minipage}{75mm}
\centering
\includegraphics[width=70mm]{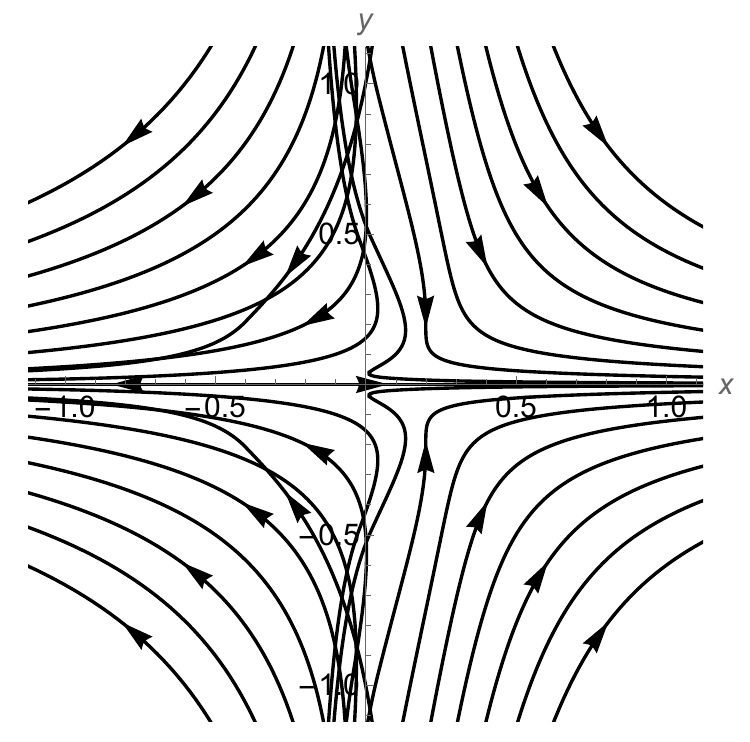}
\caption{Approximation of saddle point.}
\label{fig:1}
\end{minipage}
\end{center}
\end{figure}
\end{example}

\begin{figure}[ht]
\begin{center}
\begin{minipage}{75mm}
\centering
\includegraphics[width=70mm]{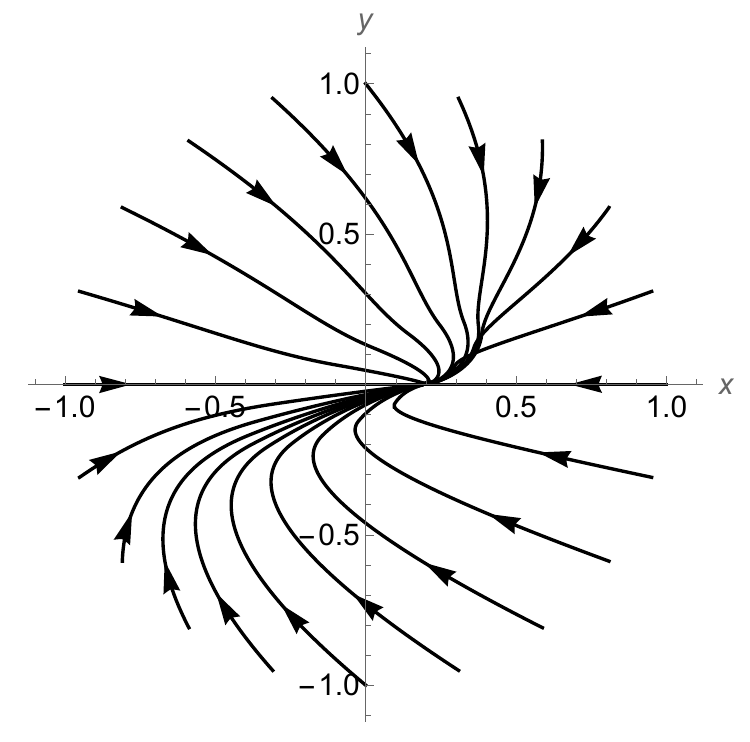}
\caption{Approximation of stable node.}
\label{fig:2}
\end{minipage}
\end{center}
\end{figure}

\begin{example}\label{ex:focus}
In \eqref{const1}, let $I=(-\infty, \infty)$, and let
\begin{equation}
 A = 
\begin{pmatrix}
-1 & 1 \\
0 & -1
\end{pmatrix}.
 \label{node}
\end{equation}
Then, by Corollary \ref{constbest2}, \eqref{const2} with \eqref{node} is Ulam stable on $I$, and the best Ulam constant is $K_{c2} = 1$. Moreover, we consider \eqref{deff} with \eqref{perturbation} and \eqref{node}. Let $\varepsilon=0.2$. Clearly the solution $\bm{\phi}(t)$ of \eqref{deff} with \eqref{perturbation} and \eqref{node} satisfies \eqref{inequ}. Since \eqref{kappa21}, \eqref{unbound21} and $\sup_{t \in I}\kappa_{21}(t) = K_{21} = K_{c2} = 1 < \infty$ hold, by Theorem \ref{UlamS2} (i), we see that there exists the unique solution $\bm{x}(t)$ of \eqref{const2} with \eqref{node} such that $\sup_{t \in I}\|\bm{\phi}(t)-\bm{x}(t)\|_{\infty} \le K_{21}\varepsilon = 0.2$. This means that there is only one orbit of \eqref{const2} with \eqref{node} in the $\varepsilon$-neighborhood of any orbit of \eqref{deff} with \eqref{perturbation} and \eqref{node}. Since the geometric classification of the origin of \eqref{const2} with \eqref{node} is a stable node (see \cite{Hale,KelPet,MicHouLiu}), a phase portrait of \eqref{deff} with \eqref{perturbation} and \eqref{node} can be proposed as an approximation of the stable node (see Figure~\ref{fig:2}). 
\end{example}

\begin{figure}[ht]
\begin{center}
\begin{minipage}{75mm}
\centering
\includegraphics[width=70mm]{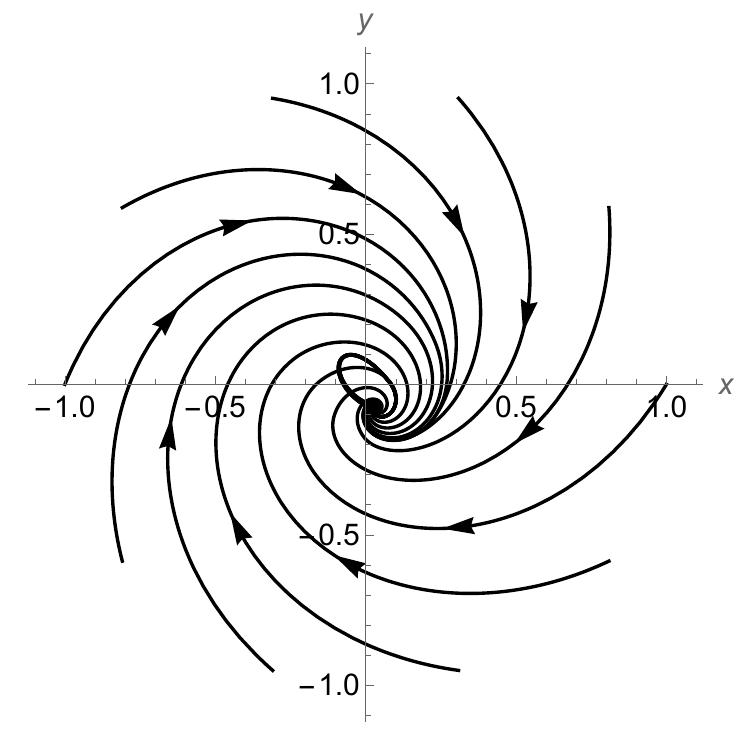}
\caption{Approximation of stable focus.}
\label{fig:3}
\end{minipage}
\end{center}
\end{figure}

\begin{example}
In \eqref{const1}, let $I=(-\infty, \infty)$, and let
\begin{equation}
 A = 
\begin{pmatrix}
-1 & 2 \\
-2 & -1
\end{pmatrix}.
 \label{focus}
\end{equation}
Then, by Corollary \ref{constbest3}, \eqref{const3} with \eqref{focus} is Ulam stable on $I$, and the best Ulam constant is $K_{c3} = 1$. Moreover, we consider \eqref{deff} with \eqref{perturbation} and \eqref{focus}. Let $\varepsilon=0.2$. Clearly the solution $\bm{\phi}(t)$ of \eqref{deff} with \eqref{perturbation} and \eqref{focus} satisfies \eqref{inequ}. Since \eqref{kappa31}, \eqref{unbound31} and $\sup_{t \in I}\kappa_{31}(t) = K_{31} = K_{c3} = 1 < \infty$ hold, by Theorem \ref{UlamS3} (i), we see that there exists the unique solution $\bm{x}(t)$ of \eqref{const3} with \eqref{focus} such that $\sup_{t \in I}\|\bm{\phi}(t)-\bm{x}(t)\|_{2} \le K_{31}\varepsilon = 0.2$. This means that there is only one orbit of \eqref{const3} with \eqref{focus} in the $\varepsilon$-neighborhood of any orbit of \eqref{deff} with \eqref{perturbation} and \eqref{focus}. Since the geometric classification of the origin of \eqref{const3} with \eqref{focus} is a stable focus (see \cite{Hale,KelPet,MicHouLiu}), a phase portrait of \eqref{deff} with \eqref{perturbation} and \eqref{focus} can be proposed as an approximation of the stable focus (see Figure~\ref{fig:3}). 
\end{example}

\section{Conclusions}
In this study, Ulam stability theorems for nonautonomous linear differential systems with generalized Jordan normal forms were established. Detailed estimates for the Ulam constants are obtained, and then the best Ulam constants are derived for all generalized Jordan normal forms. The obtained results do not require the assumption that the linear system admits an exponential dichotomy. As a result, we can get the true solution that is close to the approximate one that blows up. Analysis of Ulam stability for approximate solutions that blow up has not been done so far, so it can be said to be a new approach. In addition, the best Ulam constants for nonautonomous systems other than periodic systems have also not been presented so far. That is, this is the first study to derive the best Ulam constants for nonautonomous systems other than periodic systems. Also note that new results are obtained even in the case of constant coefficients. An example of instability is presented, as well as examples for each generalized Jordan normal form. Additionally, approximations of node, saddle, and focus are proposed.


\section*{Acknowledgments}
M. O. was supported by the Japan Society for the Promotion of Science (JSPS) KAKENHI (grant number JP20K03668) and the Research Institute for Mathematical Sciences, an International Joint Usage/Research Center located in Kyoto  University.

\end{document}